\documentclass[twocolumn]{autart}    % Enable this line and disable the 
\usepackage{graphicx} 
\usepackage{cite}
\usepackage[dvips]{epsfig}    
\usepackage[active]{srcltx}
\usepackage{algorithm}
\usepackage{algpseudocode}
\usepackage{amsmath}
\usepackage{amsfonts}
\usepackage{xcolor}
\usepackage{amssymb}
\usepackage{enumitem}
\usepackage{steinmetz}
\usepackage{tabularx}
\usepackage{hyperref}

\newcommand{\Eps}{\mathcal{E}}

\newcommand{\C}{\mathbb{C}}
\newcommand{\Z}{\mathbb{Z}}

\DeclareMathOperator*{\argmin}{arg\,min}

\begin{document}

\begin{frontmatter}
\title{Topology Learning of unknown Networked Linear Dynamical System excited by Cyclostationary inputs} 
\vspace{-0.2cm}
\author[UMN_ME]{Harish Doddi}\ead{doddi003@umn.edu},    
\author[LANL]{Deepjyoti Deka}\ead{deepjyoti@lanl.gov},         
\author[UMN_ECE]{Murti V. Salapaka}\ead{murtis@umn.edu} 
\address[UMN_ME]{Department of Mechanical Engineering, University of Minnesota Twin Cities}                \vspace{-0.2cm}          
\address[LANL]{Los Alamos National Laboratory}  
\vspace{-0.2cm}
\address[UMN_ECE]{Department of Electrical and Computer Science Engineering, University of Minnesota Twin Cities} 
\vspace{-0.2cm}        

\begin{abstract}    
Topology learning of networked dynamical systems is an important problem with implications to optimal control, decision-making over networks, cybersecurity and safety. The majority of prior work in consistent topology estimation relies on dynamical systems excited by temporally uncorrelated processes. In this article, we present a novel algorithm for guaranteed topology learning of networks that are excited by temporally (colored) cyclostationary processes, which encompasses a wide range of temporal correlation including wide-sense stationarity. Furthermore, unlike prior work, the framework applies to linear dynamic system with complex valued dependencies, and leverages group lasso regularization for effective learning of the network structure. In the second part of the article, we analyze conditions for consistent topology learning for bidirected tree networks when a subset of the network is unobserved. Here, the full topology along with unobserved nodes are recovered from observed node's time-series alone. Our theoretical contributions are validated on simulated data as well as on real-world climate data.
\end{abstract}

\end{frontmatter}

\section{INTRODUCTION}

Network representations often form an integral part in modeling the behavior of complex systems; examples include power grids \cite{deka2017structure}, thermal management systems \cite{walchli2010combined}, neuronal networks in the brain \cite{brainnetwork}, climate networks \cite{KimNorth1997}, social networks \cite{scott1988social} and the finance networks \cite{chi2010network}. Network representations constitute multiple agents which interact amongst each other constrained by a graph topology. Significant insights on the temporal and spatial evolution of the system dynamics can be gleaned using the abstractions of a network which can be employed for efficient resource management, control of assets and fault diagnosis. For instance, in epidemiology, understanding the spatial evolution of the virus spread is of interest for detecting the initial source, predicting the future impact, implementing control measures, and to limit or eradicate the virus. An important application that will be visited later in this article is about ocean currents, which are bulk movement of water from one region to another in the ocean and has a network representation. The regions are considered as nodes of the network and the current direction is represented as a network edge. Ocean currents affect the Earth's climate, navigation, rescue operations in ocean, coastal climate and economy. 

%With the rapid deployment of smart and inexpensive sensors, unveiling the interaction structure among agents from high fidelity data has become possible for many applications. 
Approaches for learning the  dependency structure of a network of dynamical systems is broadly classified into two categories: passive and active approaches. In the active learning approach, planned interventions/alterations are made to the network and the constituent agents, and the effects of the changes introduced are studied to identify network structure \cite{nabi2014network}. Passive learning, also termed non-invasive learning, learns the interaction structure without interfering with the functioning of the system \cite{friedman2008sparse, pearl2009causality}. Among others, time-series data from agents can be analyzed for exact structure learning (topology) using multivariate Wiener filtering \cite{materassi2012problem,talukdar2020physics} and directed information graphs \cite{quinn2015directed}. Note, our article focuses on learning the topology or the exact interaction graph among agents of dynamical system, and not restricted to learning the moral graph or the conditional independence graph \cite{chandrasekaran2010latent}.

Majority of the work on structure learning of dynamical systems involve systems excited by stationary processes or temporally uncorrelated or i.i.d inputs \cite{friedman2008sparse, materassi2010topological, talukdar2020physics}. However, non-stationary or temporally correlated excitation characterizes many man-made and natural phenomenon, including telecommunications \cite{communicationRef1}, seasonal weather \cite{kim1996}, biology \cite{biologyRef1}, finance \cite{stockRef1}, mechanical systems \cite{mechRef1}, and atmospheric system \cite{planetRef1}. Our prior work has studied topology learning for networks excited by wide sense stationary processes and real-valued interdependencies or dealt with real-valued time-series 
\cite{ talukdar2015reconstruction,deka2017structure,talukdar2018topology, materassi2019signal,talukdar2020physics}.

Our analysis, in this article, extends to systems where the agent's state and its interdependencies are complex-valued. An example of a complex-valued dynamical system is the power grid where nodal voltages have both real and imaginary parts \cite{deka2017structure}. More examples include, complex-valued neural networks (CVNN) where the link weights are complex-valued \cite{aizenberg2011we, hirose2013complex}, and complex-valued ordinary differential equation (CVODE) modeling \cite{yang2019complex}.
A major application of complex-valued time-series originates from spatio-temporal data such as climate data, where the fields such as velocity, temperature, pressure evolve over space and time. By taking the Fourier transform over the space dimension at various locations that are spatially distant provides complex-valued time-series.

\textbf{Contribution:} Under assumptions that hold for a large class of systems, in this article, we provide an algorithm with provable guarantees for structure identification in 
linear dynamical systems with complex-valued dependencies using nodal time-series that are modeled by cyclostationary processes. To the best of our knowledge, this is the first work to provide consistent learning in complex-valued linear dynamical systems excited by cyclostationary inputs. We also develop topology learning algorithms to partially observed bidirected {tree }networks excited by cyclostationary processes, which is of considerable importance (see \cite{ talukdar2018topology, sepehr2016inferring, sepehr2019blind}). 
We validate the theoretical contributions on the real-world climate data
as well as from the data generated from test dynamical systems. Preliminary results on some of the aspects presented in this article have appeared in the conference article \cite{talukdar2020physics, talukdar2015reconstruction, doddi2019exact}. Aside from detailed proofs of the theoretical results and extended simulation results, this article develops methods for learning with complex-valued network dependencies and under partial observability that are absent in our prior preliminary works. 
{Moreover, we provide a computational approach for our algorithms where we leverage sparse regression methods for accurate reconstruction under low sample regime.}

The remainder of the paper is organized as follows. Section \ref{sec:two} describes the topology learning problem. Section \ref{sec:three} provides the results necessary for building an algorithm to reconstruct the topology from data. Next, in Section \ref{sec:four}, an algorithm is presented for reconstructing the topology in presence of unobserved nodes for tree topologies (undirected
connected graph with no cycles). Few illustrations and applications
are provided in Section \ref{sec:five} and final conclusions in
Section \ref{sec:six}.

{\bf Notation:}
%$\boldsymbol{0}_{T \times T}$ : zero matrix of dimension $T \times T$\\ 
$A^{\prime}, A^*$ : transpose and conjugate transpose of a matrix $A$
%$A \succ 0$ : $A$ is a positive definite matrix\\
$A_{ij}\text{ or }A(i,j)\text{ or }{\boldsymbol{B}_{i}^{\prime}}A {\boldsymbol{B}_{j}}:$ $(i,j)^{th}$ element or $(i,j)^{th}$ block of size $T\times T$ in matrix $A$ (evident from the context). Here,
$\boldsymbol{B}_j = [\boldsymbol{0} \ \boldsymbol{0} \ .. \ \boldsymbol{0} \ \boldsymbol{I}_{T\times T}\ \boldsymbol{0}\ ..\ \boldsymbol{0}]^{\prime}.$ \\
%$\C$: field of complex numbers\\
$L^{2}(\Omega,\mathcal{F},\mathcal{P}): $ 
vector space of complex-valued random variables $X$ with $\mathbb{E}[X^2]< \infty,$ where $\Omega$ is a sample space, $\mathcal{F}$ is $\sigma-$algebra and $\mathcal{P}$ is a function from $\mathcal{F}$ to $[0,1]$\\
{{$\mathcal{Z}$: Z-transform\\
$x_i(k)$: $x_i$ at time $k, ~k \in \mathbb{Z}$\\
$R_x(\tau)$ : Correlation function of $x$ at lag $\tau \in \Z$\\
$\Phi_x(\omega)$: Spectral density of $x$ at frequency $\omega \in [0,2\pi)$\\
$\Re[M],~\Im [M]$: Real and imaginary part of $M$\\
$\sigma[M]:$ Largest singular value of $M$ \\
$H_{\infty}[\boldsymbol{\mathsf{W}}]: = \max_{\omega \in [0,2\pi)} \sigma[\boldsymbol{\mathsf{W}}(e^{\iota \omega})],$ for a filter $\boldsymbol{\mathsf{W}}(z)$ } }
%%%%%%%%%%%%%%%%%%%%%%%%%%%%%%%%%%%%%%%%%%%%%%%%%%%%%%%%%%%%%%%%%%%%%%%%%%%%%%%%

\section{Linear Dynamical System with cyclostationary inputs} \label{sec:two}
A Wide Sense Cyclostationary (WSCS) process of period $T$ is a random process $x(t) \in L^{2}(\Omega,\mathcal{F},\mathcal{P})$, such that  $\mathbb{E}[x(t)]=\mathbb{E}[x(t+T)],$ {$R_{xx}(s,t):= \mathbb{E}[x(s){x^*(t)}] = \mathbb{E}[x(s+T){x^*(t+T)}]$} holds for every $s,t \in \mathbb{Z}$ for a least possible  $T \in \mathbb{N}.$ Two processes $x(t)$ and $e(t),$ are {jointly wide sense cyclostationary (JWSCS)} with period $T$ if $x(t),~e(t)$ are cyclostationary with period $T$ and $R_{x,e}(s,t)=R_{x,e}(s+T,t+T)$. A wide sense stationary (WSS) process is also a WSCS process with period $T=1$. We provide a generative model for the purpose of providing theoretical guarantees. 
Suppose $\{ x_i(k)\}_{i=1}^{m}$ represents the time-series whose interaction dynamics is,
\vspace{-1.2cm}
\small
\begin{align}\label{eqn:cyclo_LDM}
{x}_i(k) &= \sum_{j=1,j\neq i}^{m}(h_{ij}*{x}_j)(k) + {e}_i(k),
\end{align}
\normalsize
where, $x_i(k)\in \C$ is the $i^{th}$ nodal time-series, $h_{ij}[k]$ is the complex-valued interdependencies, and $e_i(k)$ is an temporally correlated exogenous input (unmeasured). The transfer function $\mathsf{h}_{ij}(z):= \mathcal{Z}[h_{ij}(k)]$ with $\mathsf{h}_{ii}(z)=0, \ i,j \in \{1,\cdots,m \}$. The $\{x_j(k),e_j(k)\}_{j=1}^m$ are JWSCS of period $T$ \cite{talukdar2015reconstruction}. We refer to (\ref{eqn:cyclo_LDM}) as a networked Linear Dynamical System (LDS) and assume it is stable. We emphasize that the model (\ref{eqn:cyclo_LDM}) is applicable to systems where linearized models around an operating point {suffice to provide fidelity to reality}. Section $2$ of \cite{talukdar2020physics} provides examples of (\ref{eqn:cyclo_LDM}). However, the models considered in \cite{talukdar2015reconstruction, doddi2019exact} considered real-valued interdependencies, while it is complex-valued in our article. We will use $z$ and $\omega$ interchangeably where $z=e^{\iota \omega},\  \omega \in [0,2\pi).$

The {\it{generative graph}} $\mathcal{G} = (\mathcal{V}, \mathcal{E})$ of (\ref{eqn:cyclo_LDM}) is given by a vertex set $\mathcal{V} = \{1, ..., m\}$ and a directed edge set $\mathcal{E} = \{(i, j)|\mathsf{h}_{ij}(z) \neq 0, i,j \in \mathcal{V}\}$. Node $j$ is (i) a parent of $i$ if $(i,j)\in \mathcal{E}$, (ii) a child of $i$ if $(j,i) \in \mathcal{E},$ and (iii) a spouse of $i$ if $\exists k \in \mathcal{V}\setminus\{i,j\}$ such that $\{(k,j),(k,i)\} \in\mathcal{E}.$ The children set of $j \in \mathcal{V}$ is denoted by $\mathcal{C}_\mathcal{G}(j),$ the parent set by $\mathcal{P}_\mathcal{G}(j),$ and the spouse set by $\mathcal{K}_\mathcal{G}(j).$ The topology of $\mathcal{G}$ is an undirected graph $ \mathcal{G_T} = (\mathcal{V},  \mathcal{E_T} ),$ where $ \mathcal{E_T}:=\{(i,j)| j \in \mathcal{V}, i \in \mathcal{C}_{\mathcal{G}}(j) \cup \mathcal{P}_{\mathcal{G}}(j)\}.$ The moral graph of $\mathcal{G}$ is $\mathcal{G}_M:=(\mathcal{V}, \mathcal{E}_M)$, where $\mathcal{E}_M:=\{(i,j)| j \in \mathcal{V}, i \in \mathcal{C}_{\mathcal{G}}(j) \cup \mathcal{P}_{\mathcal{G}}(j) \cup \mathcal{K}_{\mathcal{G}}(j) \}.$ Nodes $i,j$ are in a kin relationship if $i \in \mathcal{C}_{\mathcal{G}}(j) \cup \mathcal{P}_{\mathcal{G}}(j) \cup \mathcal{K}_{\mathcal{G}}(j).$ Note that $\mathcal{E_T} \subset \mathcal{E}_M,$ and $\mathcal{E}_M \setminus \mathcal{E_T}$ is a set of strict spouses in $\mathcal{G}.$ An undirected path between $i$ and $j$ in $ \mathcal{G_T}$ is defined as a collection of distinct nodes $\{\pi_1, \cdots,\pi_k \}$ such that $\{(i,\pi_1),(\pi_1,\pi_2),\cdots, (\pi_k,j) \}\subset  \mathcal{E_T}.$ We denote $\mathcal{N}_{ \mathcal{G_T}}(i),~ \mathcal{N}_{ \mathcal{G_T}}(i,2)$ as the set of one-hop and two-hop neighbors of node $i$ in $\mathcal{G_T},$ respectively. The degree of node $i$ is $|\mathcal{N}_{ \mathcal{G_T}}(i)|.$
Two nodes $i,j \in \mathcal{V}$ are said to be {\it strict spouses} if, $j \in \mathcal{P}_G(\mathcal{C}_G(i))$ but $i\notin {\mathcal{N}}_{\mathcal{G_T}}(j).$ For a cyclostationary process $x_i(k)$ of period $T,$ the $T-$dimensional vector WSS lifted process is $X_i(k):= [x(kT)\ ...\ x(kT+T-1)]^{\prime}$. The networked LDS (\ref{eqn:cyclo_LDM}) is written in terms of vector WSS processes by lifting $\{ x_i(k)\}_{1}^m $ as below: 
\vspace{-0.8cm}
\small
\begin{align}  \label{eqn:Z_WSS}
X_j(k) &= \sum_{i=1,i \neq j}^{m}(H_{ji}*X_i)(k)+E_j(k), \text{ where, } 
\end{align}
$E_i(k) = \begin{bmatrix}e_i(kT)& e_i(kT+1)& \cdots & e_i(kT+T-1) \end{bmatrix}^{\prime},$ $\mathsf{H}_{ji}(z) = D(z^{\frac{1}{T}}) \mathsf{h}_{ji}(z^{\frac{1}{T}}),$ $D(z) = \begin{bmatrix}
    s^{\prime} & zs^{\prime} & \cdots & z^{T-1}s^{\prime} \prime
\end{bmatrix},$ and $s(z)= \begin{bmatrix} z^0 & z^{-1} &\cdots &z^{-(T-1)} \end{bmatrix}.$
\normalsize
The $z-$transform of (\ref{eqn:Z_WSS}) is  $\mathsf{X}_{j}(z) = \sum_{i=1,i \neq j}^{m}\mathsf{H}_{ji}(z)\mathsf{X}_i(z) + \mathsf{E}_j(z).$ The proof of (\ref{eqn:Z_WSS}) is given in Lemma \ref{lem:lemma4} in Appendix \ref{app:lemma4_proof}. By stacking the vector processes of (\ref{eqn:Z_WSS}) in $z-$ domain, we get:
\vspace{-0.6cm}
\small
\begin{align} \label{vector_LDM_gen}
z-\text{domain: }&\mathsf{X}(z) = \mathbb{H}(z)\mathsf{X}(z)+ \mathsf{E}(z), 
\end{align}
\normalsize
where, $~ \forall i,j \in \mathcal{V},$ the $(i,j)^{th}$ block of $\mathbb{H}(z)$ is $\mathsf{H}_{ij}(z),$ $\mathsf{X} = [\mathsf{X}_1, \cdots, \mathsf{X}_m]^{\prime},$
$\mathsf{E} = [\mathsf{E}_1, \cdots, \mathsf{E}_m]^{\prime}.$ Since,  cyclostationary processes $\{e_i(k) \}_{i=1}^m$ are mutually uncorrelated, it follows that $\Phi_E$ is block diagonal. We assume  (\ref{vector_LDM_gen}) is \emph{well-posed and topologically detectable}, that is, $\Phi_{E}^{-1}(\omega) \succ 0$ and $\mathbb{I}_{mT \times mT}-\mathbb{H}(\omega)$ is invertible almost surely for $\omega \in[0,2\pi),~i=1,...,m$. It follows from (\ref{eqn:Z_WSS}), a directed edge from $i$ to node $j$ exists in $\mathcal{G}$ if and only if $\mathsf{H}_{ji} \neq \mathbf{0}$ in (\ref{eqn:Z_WSS}). Thus, LDG for the lifted vector WSS processes is identical to the LDG of (\ref{eqn:cyclo_LDM}) \cite{talukdar2015reconstruction}; we provide an algorithm to learn the $\mathcal{G_T}$ using the lifted vector WSS processes.

%\textcolor{red}{given the first problem is for general LDMs and second is for tree , should you say that upfront and also not say that the determine $\mathcal{G}_\mathcal{T}$ for the first problem (it isnt a tree)}. 

\section{Learning $\mathcal{G_T}$ for a fully observed network} \label{sec:three}
{\it Given $\{ x_i(k) \in \mathbb{C}, k\in \mathbb{Z}\}, $ for all $i \in \mathcal{V}$ of networked LDS $(\mathbb{H}(z),{E})$ described by (\ref{vector_LDM_gen}) which is well-posed and topologically detectable, determine $ \mathcal{G_T}(\mathcal{V},  \mathcal{E_T}).$ Directed loops in $\mathcal{G}(\mathcal{V}, \mathcal{E})$
are admissible.} We provide the topology learning results in terms of the properties of inverse power spectral density of $\mathsf{X}$ in (\ref{vector_LDM_gen}) and the multivariate Wiener filtering.
The learning of $\mathcal{E_T}$ involves the reconstructing $\mathcal{E}_M$ and then eliminating the strict spouse edges ($\mathcal{E}_M \setminus \mathcal{E_T}$) from $\mathcal{E}_M.$ The reconstruction of $\mathcal{G}_M(\mathcal{V},\mathcal{E}_M)$ follows from the support structure of $\Phi^{-1}_X,$ as 
shown in Section III of  \cite{doddi2019exact}. From ($13$) in \cite{doddi2019exact}, it follows that,
\vspace{-0.5cm}
\small
\begin{align}\label{PhiInv_expression}
\hspace{-0.4cm}{\boldsymbol{B}_{j}^{\prime}}{{\Phi}_{X}^{-1}}{\boldsymbol{B}_{i}}=-{{\Phi}_{E_j}^{-1}}{\mathsf{H}_{ji}}-{{\mathsf{H}^{*}_{ij}}}{{\Phi}_{E_i}^{-1}}+{\sum_{k=1}^{m}}{(\mathsf{H}_{kj})^{*}}{{\Phi}_{E_k}^{-1}}(\mathsf{H}_{ki}).
\end{align}
\normalsize
If ${\boldsymbol{B}_{j}^{\prime}}{{\Phi}_{X}^{-1}}{\boldsymbol{B}_{i}} \neq \boldsymbol{0},$ then $i \in \mathcal{C}_{\mathcal{G}}(j) \cup \mathcal{P}_{\mathcal{G}}(j) \cup \mathcal{K}_{\mathcal{G}}(j)$ (kin relationship). However, the converse is not guaranteed and such cases are pathological (see \cite{talukdar2015reconstruction, doddi2019exact}).
Consider an edge set $\bar{\mathcal{E}}_M := \{(i,j)|{\boldsymbol{B}_{j}^{\prime}}{{\Phi}_{X}^{-1}}{\boldsymbol{B}_{i}} \neq \boldsymbol{0} \}.$ Thus, $\bar{\mathcal{E}}_M = \mathcal{E}_M.$ However, the number of spurious edges in $\mathcal{E}_M$, that is $\mathcal{E}_M \setminus  \mathcal{E_T},$ can be substantial. 
%We provide a necessary and sufficient condition to identify the strict spouse edges in $\mathcal{G}(\mathcal{V}, \mathcal{E})$ to eliminate them from the reconstructed $\bar{\mathcal{E}}_M$ to recover the topology $\mathcal{E_T}$. 
\subsection{Identifying strict spouses in $\mathcal{G}(\mathcal{V}, \mathcal{E})$ }
{
To provide a sufficient condition for identifying strict spouses in $\mathcal{G},$ we make the following assumption on the number of strict spouses in the generative graph. 

\begin{assum}\label{ass:atmost1}
In the generative graph $\mathcal{G}$ of the networked LDS of (\ref{eqn:cyclo_LDM}), for any two strict spouses $i$ and $j$, the set of all common children $\{k|k \in \mathcal{C}_G(j)\cap \mathcal{C}_G(i)\}$ has a cardinality of at most $1.$
\end{assum}

Assumption \ref{ass:atmost1} is not necessary when all the interdependencies in (\ref{eqn:cyclo_LDM}) are strictly real-valued, a scenario presented in \cite{doddi2019exact}. The above assumption is satisfied by a large class of directed graphs, such as with (i) tree topologies (undirected connected graph with no cycles), (ii) loopy topologies, with every loop has a size greater than four. Such networks are common in infrastructure networks such as power and gas grids, as the practical loop size is much greater than four. If Assumption \ref{ass:atmost1} is not satisfied, then the reconstructed topology would contain spurious edges that are restricted to within two-hop distance from a node in $\mathcal{G}.$

\begin{thm}\label{thm:3.2}
Consider a well-posed and topologically detectable networked LDS described by (\ref{eqn:cyclo_LDM}), with its equivalent representation $(\mathbb{H}(z),{E})$ as in (\ref{vector_LDM_gen}), its associated LDG $\mathcal{G},$ topology $\mathcal{G_T}$ and satisfying Assumption \ref{ass:atmost1}. Let $i,j \in \mathcal{V},$ are strict spouses in $\mathcal{G}(\mathcal{V}, \mathcal{E})$. Suppose $\Eps_{ji} \in \C^{T \times 1}$ be the eigenvalues of ${\boldsymbol{B}_{j}^{\prime}}{{\Phi}_{X}^{-1}(\omega)}{\boldsymbol{B}_{i}}.$ Then, $\phase{\Eps_{ji}(l)}$ is a constant, independent of $\omega,$ for $l\in \{1,2,\cdots,T \}$ and for all $\omega \in [0, 2 \pi)$.
\end{thm}
\begin{pf}
Consider $j \in \mathcal{P}_G(\mathcal{C}_G(i))$ and $i\notin \mathcal{N}_{\mathcal{G_T}}(j)$, then ${\mathsf{H}_{ji}(\omega)}= {\mathsf{H}_{ij}(\omega)}= \boldsymbol{0}_{T \times T}$ for all $\omega \in [0,2\pi).$ Using (\ref{PhiInv_expression}),

\small
\begin{align} 
&{\boldsymbol{B}_{j}^{\prime}}{{\Phi}_{X}^{-1}(\omega)}{\boldsymbol{B}_{i}}=-{{\Phi}_{E_j}^{-1}(\omega)}{\mathsf{H}_{ji}(\omega)}-({\mathsf{H}_{ij}(\omega)})^*{{\Phi}_{E_i}^{-1}(\omega)}+\nonumber \\
&\ \ \ \ \ \ \ \ \ \ \ \ \ \ \ {\sum_{k \in \mathcal{C}_G(j)\cap \mathcal{C}_G(i)}}(\mathsf{H}_{kj}(\omega))^*{{\Phi}_{E_k}^{-1}(z)}\mathsf{H}_{ki}(\omega)\nonumber \\
&= {\sum_{k \in \mathcal{C}_G(j)\cap \mathcal{C}_G(i)}}{(\mathsf{H}_{kj}(\omega))^{*}}{{\Phi}_{E_k}^{-1}(\omega)}\mathsf{H}_{ki}(\omega) \nonumber \\
&={\sum_{k \in \mathcal{C}_G(j)\cap \mathcal{C}_G(i)}}[D(\frac{\omega}{T})\mathsf{h}_{kj}(\frac{\omega}{T})]^*{{\Phi}_{E_k}^{-1}(\omega)}D(\frac{\omega}{T})\mathsf{h}_{ki}(\frac{\omega}{T}).\nonumber \\
%&(\because \text{ using (\ref{eqn:H_h})})\nonumber \\
&={\sum_{k \in \mathcal{C}_G(j)\cap \mathcal{C}_G(i)}}[h^*_{kj}(\frac{\omega}{T})h_{ki}(\frac{\omega}{T})][D^*(\frac{\omega}{T}){\Phi}_{E_k}^{-1}(\omega)D(\frac{\omega}{T})] \nonumber\\
&=[h^*_{k_1j}(\frac{\omega}{T})h_{k_1i}(\frac{\omega}{T})][D^*(\frac{\omega}{T}){\Phi}_{E_{k_1}}^{-1}(\omega)D(\frac{\omega}{T})] \nonumber\\
&(\because \text{using Assumption \ref{ass:atmost1}, where $k_1 = \mathcal{C}_G(j)\cap \mathcal{C}_G(i)$} )\nonumber \\
&=|h_{k_1j}(\frac{\omega}{T})|^2[\frac{h_{k_1i}(\frac{\omega}{T})}{h_{k_1j}(\frac{\omega}{T})}] \times \text{ Positive Definite Hermitian. }
\end{align}
\normalsize

{Define $a_{ik_1j} = \frac{h_{k_1i}(\frac{\omega}{T})}{h_{k_1j}(\frac{\omega}{T})}.$ The eigenvalues of a positive definite hermitian are positive real-valued. Suppose $\Eps_{ji} \in \C^{T \times 1}$ denotes the eigenvalues of ${\boldsymbol{B}_{j}^{\prime}}{{\Phi}_{X}^{-1}(\omega)}{\boldsymbol{B}_{i}}.$ Then, $\phase{\Eps_{ji}(l)}=\phase{a_{ik_1j}},$ a constant independent of $\omega \in [0,2\pi)$ for $l\in \{1,2,\cdots,T \}.$}
    
\end{pf}

The Theorem \ref{thm:3.3} from Appendix \ref{app:thm3.3_proof} shows that the converse of Theorem \ref{thm:3.2} holds, except for pathological cases. Hence, Theorem \ref{thm:3.2} is used as a necessary and sufficient condition to prune out spurious edges from $\bar{\mathcal{E}}_M$ or $\mathcal{E}_M.$

\subsection{Relaxing Assumption \ref{ass:atmost1} - applicable for any generative graph}

For identifying strict spouses in generative graphs when multiple common children between strict spouses are allowed, that is, Assumption \ref{ass:atmost1} is relaxed, we make the following assumption to develop a condition for identifying the strict spouses.  

\begin{assum} \label{ass:transferFunction_property}
For $\mathsf{h}_{ki}(\omega), \mathsf{h}_{kj}(\omega) \neq 0,$ then
$\phase{\mathsf{h}_{kj}(\omega)} - \phase{\mathsf{h}_{ki}(\omega)}$ is a constant for all $\omega \in [0,2\pi)$.
\end{assum}

Define
$a_{jki} := \frac{\mathsf{h}_{kj}}{\mathsf{h}_{ki}},$ for $\mathsf{h}_{ki}(\omega), \mathsf{h}_{kj}(\omega) \neq 0.$  
From Assumption \ref{ass:transferFunction_property}, it follows that for strict spouses $i,j \in \mathcal{V},$ the 
$\phase{a_{jki}}$ is a constant for all $\omega \in [0,2\pi)$. Assumption \ref{ass:transferFunction_property} is  satisfied by a large class of dynamical systems, including but not restricted to the generative models considered in \cite{doddi2019exact, doddi2022efficient,  materassi2012problem, talukdar2017exact, talukdar2020physics,  talukdar2015reconstruction}. From the term ${\sum_{k=1}^{m}}{(\mathsf{H}_{kj})^{*}}{{\Phi}_{E_k}^{-1}}(\mathsf{H}_{ki})$ in (\ref{PhiInv_expression}), it follows that Assumption \ref{ass:transferFunction_property} holds trivially for a network with $\{ a_{ikj} \in \mathbb{R}\}, $ a scenario considered in \cite{doddi2019exact, doddi2022efficient,  materassi2012problem, talukdar2017exact, talukdar2020physics,  talukdar2015reconstruction}. An illustration is shown in Appendix \ref{app:illustration}. However, for complex-valued interdependencies in the generative model, Assumption \ref{ass:transferFunction_property} is a sufficient condition for identifying strict spouse edges without restricting the structure of the generative graph. Moreover, if the network comprises of WSS processes (time-period $T=1$), the ${{\Phi}_{X}^{-1}(\omega)}(j,i)$ is a scalar quantity (not a matrix). Then, for strict spouses $i,j \in \mathcal{V},$  $\phase{{{\Phi}_{X}^{-1}(\omega)}(j,i)}$ is a constant $\theta$ for all $\omega \in [0, 2 \pi),$ a result presented in \cite{doddi2022efficient,  talukdar2017exact, talukdar2020physics}. 

\begin{thm}\label{thm:3.4}
Consider a well-posed and topologically detectable networked LDS described by (\ref{eqn:cyclo_LDM}),  its representation $(\mathbb{H}(z),{E})$ as in (\ref{vector_LDM_gen}), associated LDG $\mathcal{G}$ and topology $\mathcal{G_T},$ while satisfying Assumption \ref{ass:transferFunction_property}. Suppose, for {$i,j \in \mathcal{V}, \phase{a_{ik_1j}}=\phase{a_{ik_2j}}$} holds for all $k_1,k_2 \in \mathcal{C}_G(i) \cap \mathcal{C}_G(j), k_1 \neq k_2.$ Let $\Eps_{ji} \in \C^{T \times 1}$ be the eigenvalues of ${\boldsymbol{B}_{j}^{\prime}}{{\Phi}_{X}^{-1}(\omega)}{\boldsymbol{B}_{i}}$, where $X(k)$ is the output of the networked LDS (\ref{vector_LDM_gen}). Then, $\phase{\Eps_{ji}(l)}$ is a constant $c_{ji}$ for all $l\in \{1,2,\cdots,T \},$ $\omega \in [0, 2 \pi)$, if and only if $j \in \mathcal{P}_G(\mathcal{C}_G(i))$ but $i\notin \mathcal{N}_{\mathcal{G_T}}(j)$, that is, $i,j$ are strict spouses.
\end{thm}

\begin{pf}
Suppose {$a_{ikj} = r_ke^{\iota \theta_k}$} for all $k \in \mathcal{C}_G(j)\cap \mathcal{C}_G(i).$ Under the stated assumption, we have $\phase{a_{ik_1j}}=\phase{a_{ik_2j}}$ for any distinct $k_1,k_2\in \{k|k \in \mathcal{C}_G(j)\cap \mathcal{C}_G(i)\}$. Then, $\theta_k$ is independent of $k.$ Let $\theta_k = \theta$ for all $k \in \mathcal{C}_G(j)\cap \mathcal{C}_G(i).$
For strict spouses $i,j$, we have,
 \small
 \begin{align*}
 &{\boldsymbol{B}_{j}^{\prime}}{{\Phi}_{X}^{-1}(\omega)}{\boldsymbol{B}_{i}}\\
&={\sum_{k \in \mathcal{C}_G(j)\cap \mathcal{C}_G(i)}}[a_{ikj}|h_{kj}(\frac{\omega}{T})|^2][D^*(\frac{\omega}{T}){\Phi}_{E_k}^{-1}(\omega)D(\frac{\omega}{T})] \nonumber\\
 &= {\sum_{k \in \mathcal{C}_G(j)\cap \mathcal{C}_G(i)}}{r_ke^{\iota \theta_k}}[D^*(\frac{\omega}{T}){\Phi}_{E_k}^{-1}(\omega)D(\frac{\omega}{T})]\\
 &= e^{\iota \theta}{\sum_{k \in \mathcal{C}_G(j)\cap \mathcal{C}_G(i)}}r_k[D^*(\frac{\omega}{T}){\Phi}_{E_k}^{-1}(\omega)D(\frac{\omega}{T})] \\
 & = e^{\iota \theta}\times \text{ Positive Definite Hermitian. }\\
 &\text{Hence, $\phase{\Eps_{ji}(l)}$ is $\theta,\ \forall l\in \{1,2,\cdots,T \}.$}
 \end{align*}
 \normalsize
 
Here, $c_{ji}:= \theta.$ The converse of the Theorem \ref{thm:3.4} holds, except for pathological cases. The proof is similar to Theorem \ref{thm:3.3}. Therefore, the consequences of Theorem \ref{thm:3.2} hold even though Assumption \ref{ass:atmost1} is violated.
    
\end{pf}

% \small
% \begin{algorithm}
% \caption{Learning Algorithm for reconstructing the topology of LDG with cyclostationary inputs}
% \textbf{Input:} Nodal time-series $x_i(k)$ for each node $i \in \{1, 2,... m\}$ which is WSCS. Thresholds $\rho,\tau$. Frequency points $\Omega$. \\
% \textbf{Output:} Reconstructed Topology ($\mathcal{V},\bar{\mathcal {E_T}}$) \\
% \begin{algorithmic}[1]
% \State Compute the periodogram of $x_i(k)$ to determine the period $T_i$. Determine $T= LCM\{T_1,\cdots, T_m \}.$ Lift each $x_i(k)$ to vector process $X_i(k)$ of size $T.$ Define $X(k) = [X_1(k),\cdots X_m(k)]^{\prime}$
% \State Edge set $\bar{\mathcal{E}}_M \gets \{\}$ 
% \ForAll{$l,p \in \{1,2,...,m\}, l\neq p$}
% \If{$H_{\infty}[ \boldsymbol{\mathsf{W}}_{lp} ] \neq 0$}
% \State $\bar{\mathcal{E}}_M \gets \bar{\mathcal{E}}_M \cup \{(l,p)\}$
% \EndIf
% \EndFor\label{step1_b}
% \State Edge set $\bar{\mathcal{E}}_T \gets \bar{\mathcal{E}}_M$ 
% \ForAll{$(p,l) \in \bar{\mathcal{E}}_M$}
% \State Compute $\{{\Eps_{pl}(t)\}}^{\prime}_{t=1}=eig\{\boldsymbol{\mathsf{W}}_{pl}(\omega)\}$
% \If{$\phase{\Eps_{pl}(t)}$ is constant $\forall \omega \in [0,2\pi), \forall t$}
% \State $\bar{\mathcal{E}}_T \gets \bar{\mathcal{E}}_T - \{(p,l)\}$
% \EndIf
% \EndFor 
% \end{algorithmic}
% \end{algorithm}
% \normalsize

\subsection{Connection between multivariate Wiener filter and Power spectral density}
The connection between inverse power spectral density and the multivariate Wiener filtering for topology learning of cyclostationary processes has been explored in \cite{talukdar2015reconstruction,doddi2019exact}.
Using ($12$) from \cite{doddi2019exact}, we have,
\vspace{-0.7cm}
\small
\begin{align}\label{eqn:Phi_wji}
 &{\boldsymbol{B}_{j}^{\prime}}{{\Phi}_{X}^{-1}}(\omega){\boldsymbol{B}_{i}}=-\Phi_{\Eps_j}^{-1}\boldsymbol{\mathsf{W}}_{ji}(\omega),
 % &\text{where } \Eps_j(k)= X_j(k)- \sum_{i=1,i\neq j}^{m}(\boldsymbol{W}_{ji}*X_i)(k).\nonumber
\end{align}
\normalsize
where $ \Eps_j(k)= X_j(k)- \sum_{i=1,i\neq j}^{m}(\boldsymbol{W}_{ji}*X_i)(k).$ Using (\ref{eqn:Phi_wji}) in (\ref{PhiInv_expression}), the moral graph $\mathcal{G}_M$ can be reconstructed using the condition $\boldsymbol{\mathsf{W}}_{ji}(\omega) \neq \boldsymbol{0}$ instead of ${\boldsymbol{B}_{j}^{\prime}}{{\Phi}_{X}^{-1}}(\omega){\boldsymbol{B}_{i}} \neq \boldsymbol{0}.$ Note that $\Phi_{\Eps_j}^{-1}$ in (\ref{eqn:Phi_wji}) is positive definite and its eigenvalues are real. From (\ref{eqn:Phi_wji}), the phase of eigenvalues of $\boldsymbol{\mathsf{W}}_{ji}(\omega)$ is equal to phase of eigenvalues of ${\boldsymbol{B}_{j}^{\prime}}{{\Phi}_{X}^{-1}}(\omega){\boldsymbol{B}_{i}}+\pi.$ 

%For Theorem \ref{thm:3.2}, we compute the eigenvalues of $\boldsymbol{\mathsf{W}}_{ji}(\omega),$ for all $(i,j) \in \mathcal{E}_M, $ and check if the phase of the eigenvalues remain constant for all $\omega \in [0,2\pi).$ 

\subsection{Regularized Wiener Filter Estimate:} \label{subsec:RWFE}

Wiener filter estimation from the data is a necessary step for both fully observability and partially observability of the network. 
{\it The theoretical results provided in this article are based on the asymptotic Wiener filter estimator}. Sparse regression techniques such as Group Lasso regularization \cite{bach2008consistency} is used in the Wiener filter computation for topology reconstruction in the high dimensional setting. So, we propose a block-sparsity Lasso regularized problem to promote the sparsity among the $T \times T$ blocks in the Wiener filter estimate matrix. The regularized Wiener filter can be estimated both in time-domain as well as in frequency domain.

\textbf{Time-domain regression:}
The following time-domain optimization is proposed for a nodal time-series trajectory $\{X_j(k) \in \mathbb{C}\}$:
\vspace{-0.8cm}
\small
\begin{align}\label{eqn:vec_opt_sparse}
	&{\min_{{\{\boldsymbol{W}_{ji}^L\}}_{ L = -L_*}^{L = L_*}}}~~ \mathbb{E}\big[\| X_j(k)-\sum_{i=1,i\neq j}^{m}  \sum_{L=-L_*}^{L=L_*}\boldsymbol{W}_{ji}^L X_i(k-L) \|_2^2 \big] \nonumber \\ & ~~~~~~~~~~~~~~~~~~~~~~~ +\gamma \sum_{i=1,i \neq j}^{m}\sqrt{\sum_{L=-L_*}^{L_*} \| \boldsymbol{W}_{ji}^L \|_F^2}, 
\end{align}
\normalsize
where, $\gamma$ is the regularization parameter, $L_*$ is the maximum lag, and $\|.\|_F$ is the Frobenius norm. Both $\gamma$ and $L_*$ are predetermined constants tuned based on historical data. From the optimal solution of (\ref{eqn:vec_opt_sparse}), we construct 
$\boldsymbol{\hat{\mathsf{W}}}_{ji}(\omega)= \sum_{L= -L_*}^{L_*} \boldsymbol{W}_{ji}^{L}e^{-\iota \omega L},$ and use $\boldsymbol{\hat{\mathsf{W}}}_{ji}(\omega)$ as a surrogate to $\boldsymbol{\mathsf{W}}_{ji}(\omega).$ 
In our preliminary work \cite{doddi2019exact}, no regularization is employed ($\gamma =0$), and thus (\ref{eqn:vec_opt_sparse}) is solved as $T$ independent scalar optimization problems. We remark that it requires relatively small number of samples when (\ref{eqn:vec_opt_sparse}) is solved with regularization ($\gamma \neq 0$). {\it Note that (\ref{eqn:vec_opt_sparse}) cannot be solved as $T$ independent scalar optimization problems due to the coupling present in the regularizer term.} In Section \ref{sec:five}, we demonstrate the performance of reconstructing topology based on (\ref{eqn:vec_opt_sparse}) for an ocean current network. 

\textbf{Frequency-domain regression:}
Consider a nodal time-series trajectory $\{ X_j(k) \in \C\},$ $k \in \{1, \cdots, n\times N \}$ containing $n \times N$ timestamps for a node $j \in \mathcal{V}$. We split the nodal trajectory $\{ X_j(k) \}_{k=1}^{n \times N}$ into $n$ sample trajectories each of length $N.$ Let the 
$r^{th}$ sample trajectories is denoted by $\{X_j^r(k) \in \C \}_{k = 1}^N,$ $r \in \{1, \cdots, n \}.$ For the $r^{th}$ sample trajectory, define the Discrete Fourier Transform (DFT)\footnote{computed at frequency $\omega = \frac{2\pi l}{N},~ l \in \{0,\cdots, N-1 \}$ unless explicitly mentioned} as
\vspace{-0.6cm}
\small
\begin{align}\label{def:X_i}
 &{\tilde{X}}^r_j = \frac{1}{\sqrt{N}}\sum_{k=1}^{N} X_j^r(k) e^{-\iota \omega k},  {\tilde{X}}^r_{\overline{j}} = \frac{1}{\sqrt{N}}\sum_{k=1}^{N} X_{\overline{j}}^r(k) e^{-\iota \omega k},
\end{align}
\normalsize
where, $ X_{\overline{j}}^r =[(X_1^r)^{\prime}, \cdots,(X_{j-1}^r)^{\prime},(X_{j+1}^r)^{\prime}, \cdots,(X_{m}^r)^{\prime}]^{\prime},$ $r \in\{1, \cdots,n \}.$ Construct $\mathcal{Y} = [{\tilde{X}}_j^1, \cdots, {\tilde{X}}_j^n ]^{\prime} \in \C^{n \times T}$ and $\mathcal{X} = [{\tilde{X}}_{\overline{j}}^1, \cdots, {\tilde{X}}_{\overline{j}}^n]^{\prime} \in \C^{n \times (m-1)T}$ respectively. 
% In the asymptotic limit, the true Wiener filter for a node $j \in \mathcal{V}$ is obtained by solving the below optimization problem:
% \vspace{-0.5cm}
% \small
% \begin{align} \label{eqn:RWF}
%  \boldsymbol{\mathsf{W}}_j(\omega) = \argmin_{\beta \in \mathbb{C}^{(m-1)T \times T}} \lim_{n,N \to \infty} \frac{1}{2n}\| \mathcal{Y} - \mathcal{X}\beta \|_2^2.
% \end{align}
% \normalsize
For finite $N$, we propose a regularized version of Wiener filter estimate:
\vspace{-0.7cm}
\small
\begin{align} \label{eqn:RWFE} \hat{\boldsymbol{\mathsf{W}}}_j(\omega, \lambda) = &\argmin_{\beta= [(\tilde{\beta}_1)^{\prime}, \cdots, (\tilde{\beta}_{m-1})^{\prime}]^{\prime}, 
 \tilde{\beta}_l \in \mathbb{C}^{T \times T}} \frac{1}{2n}\|\mathcal{Y} - \mathcal{X}\beta \|_2^2  \nonumber \\
 &~~~~~~~~+\lambda \sum_{l=1}^{m-1} \| \tilde{\beta}_l \|_F,
\end{align}
\normalsize
where, $\lambda >0$ is the regularization parameter, $\|. \|_F$ is the Frobenius norm of a matrix. The entries of the matrix $(\tilde{\beta}_l)^{\prime}$ belong to a group and the regularization term in (\ref{eqn:RWFE}) promotes block sparsity in $\hat{\boldsymbol{\mathsf{W}}}_j.$ We demonstrate the applicability of (\ref{eqn:RWFE}) in reconstructing the network of ocean currents from sea surface temperatures (SST).
%\subsection{Summary of Topology Learning under Full Observability:}
%We propose Algorithm $1$ for reconstructing $ \mathcal{G_T}$ from time-series $\{x_i(k) \}_{i=1}^m.$ The Algorithm $1$ shown in Appendix \ref{app:alg1}. It first estimates time-period $T$ and lifts each times series to vector WSS process (Steps $1-2$). Then, multivariate Wiener filter is computed from the lifted processes (Steps $3-5$). (\ref{PhiInv_expression}) is then used to reconstruct the moral graph (Steps $6-11$).
%Further, if Assumptions \ref{ass:transferFunction_property} and \ref{ass:atmost1} hold, then the spurious edges from the reconstructed moral graph are identified and eliminated using Theorems \ref{thm:3.2}, \ref{thm:3.3} and Remark \ref{rem:pathology1}. The true topology is thus reconstructed.

%We discussed the topology learning when the measurements of all the nodes in the network are available. However, often in many scenarios some of the nodes might be unobserved or latent. Then, topology learning also involves extracting the latent nodes and their neighbors. We discuss the topology learning with latent nodes in the next section.

{
\section{Topology Learning of Cyclostationary time-series with Latent nodes} \label{sec:four}}
We focus on exact topology learning when only a subset of nodes in the network are observed (partial observability), $\mathcal{G}(\mathcal{V},\mathcal{E})$ is bidirected and $\mathcal{G_T}$ is tree, with complex-valued interdependencies. Suppose $\mathcal{V}_o:=\{1,\cdots, m \}$ be the set of observed nodes and $\mathcal{V}_h:=\{m+1,\cdots, n \}$ be the set of latent nodes. 
Given $\{ x_i(k), k\in \mathbb{Z}\}$ for all $i \in \mathcal{V}_o \subset \mathcal{V}$ of networked LDS $(\mathbb{H}(z),{E})$ described by (\ref{vector_LDM_gen}) which is well-posed and topologically detectable, we determine $ \mathcal{G_T}(\mathcal{V},  \mathcal{E_T}).$ WLOG, suppose $\{x_i(k) \}_{i=1}^m$ are observed processes and $\{x_i(k) \}_{i=m+1}^n$ are the latent processes.
The exogenous inputs $\{e_i(k) \}_{i=1}^{n}$ are mutually uncorrelated WSCS of period $T:= LCM\{T_1,\cdots,T_n \}.$ 

\begin{assum}\label{ass:timeperiod}
The time-period $T_l,$ for all $l \in \{m+1, \cdots,n \},$ can be written as $\frac{T}{n}$ for some $n \in \mathbb{N}.$
\end{assum}

By Assumption \ref{ass:timeperiod}, we have {$T = LCM\{T_1, \cdots, T_m, \cdots, T_n \}$ } $= LCM\{T_1, \cdots, T_m \}.$ Then, the collection $\{x_i(k),e_i(k) \}_{i=1}^{n}$ are jointly cyclostationary with period $T,$ which can be solely computed from observed nodal time-series. Assumption \ref{ass:timeperiod} is not required in Section \ref{sec:three} because all the nodes are observable. All processes $\{x_i(k), e_i(k)\}_{i=1}^n$ are lifted to a $T-$dimensional vector WSS process $\{X_i(k),E_i(k)\}_{i=1}^n.$ The network dynamics of $\{X_i(k)\}_{i=1}^n$ is,
\vspace{-0.7cm}
\small
\begin{align}\label{eqn:temp}
\left[ \begin{array}{c} X_o(k) \\ X_h(k) \end{array} \right] &=  \begin{bmatrix} \mathsf{H}_{oo} & \mathsf{H}_{oh} \\ \mathsf{H}_{ho} & \mathsf{H}_{hh} \end{bmatrix}*\left[ \begin{array}{c} X_o \\ X_h \end{array} \right] (k) + \begin{bmatrix} E_o(k) \\ E_h(k)\end{bmatrix}, 
\end{align}
\normalsize
\vspace{-0.7cm}
\small where, $X_o=\begin{bmatrix} X_1, \cdots,  X_m \end{bmatrix}^{\prime},~X_h=\begin{bmatrix} X_{m+1},  \cdots, X_n \end{bmatrix}^{\prime},$ 
$E_o= \begin{bmatrix} E_1, \cdots, E_m \end{bmatrix}^{\prime}, E_h= \begin{bmatrix} E_{m+1}, \cdots, E_n \end{bmatrix}^{\prime}.$ \normalsize 
The dynamics (\ref{eqn:temp}) is written in a compact form as,
\vspace{-0.65cm}
\small 
\begin{align} \label{eqn:latent_LDM}
X(k) = \mathbb{H}(z)*X(k)+ E(k).    
\end{align}
\normalsize
Here, $\mathsf{H}_{oo}$ is $mT \times mT$ transfer matrix representing the interactions from the observed nodes to all the observed nodes. $\mathsf{H}_{oh}$ is of size $mT \times (n-m)T$ represents the interactions from the latent nodes to all the observed nodes. $\mathsf{H}_{ho}$ is of size $(n-m)T \times mT,$ and it represents the interaction term from the observed nodes to all the latent nodes. Finally, $\mathsf{H}_{hh}$ is of size $(n-m)T \times (n-m)T$ and it represents the interactions from the latent nodes to all the latent nodes. Define $\mathcal{E}= \{(i,j)| i,j \in \mathcal{V},\ \mathbb{H}(i,j) \neq \boldsymbol{0} \}.$ We denote the generative graph associated with (\ref{eqn:latent_LDM}) by $\mathcal{G}(\mathcal{V},\mathcal{E})$ and its topology by $ \mathcal{G_T}(\mathcal{V}, \mathcal{E_T}).$ 
Further, the restriction of $\mathcal{G}(\mathcal{V},\mathcal{E})$ to the observed nodes is denoted by $\mathcal{G}_o(\mathcal{V}_o,\mathcal{E}_o),$ where $\mathcal{E}_o:=\{(i,j)|i,j \in \mathcal{V}_o, (i,j) \in \mathcal{E} \}.$ The remaining directed edges in $\mathcal{E}$ are associated with hidden nodes, denoted by $\mathcal{E}_h:= \mathcal{E}\setminus \mathcal{E}_o.$ The topology $ \mathcal{G_T}$ restricted to $\mathcal{V}_o$ is $\mathcal{G}_{T_o}(\mathcal{V}_o,\mathcal{E}_{T_o}),$ where the undirected edge set $\mathcal{E}_{T_o}:=\{(i,j)|i,j \in \mathcal{V}_o, (i,j)\text{ or }(j,i) \in \mathcal{E}_o \}.$ Define $\mathcal{E}_{T_h}:=  \mathcal{E_T} \setminus \mathcal{E}_{T_o}.$ 

{\it Given the nodal observed time-series $\{x_i(k)\}_{i=1}^m,$ we determine 
the topology $ \mathcal{G_T}(\mathcal{V},  \mathcal{E_T}).$ } Our framework for learning $ \mathcal{G_T}$ comprises reconstructing (a) the edges between observed nodes $\mathcal{G}_{T_o}(\mathcal{V}_o,\mathcal{E}_{T_o})$ (b) find hidden nodes ($\mathcal{V}_h$) and (c) the connections to hidden nodes and the observed nodes ($\mathcal{E}_{T_h}$). From (\ref{eqn:latent_LDM}), it follows that,
\vspace{-0.6cm}
\small
\begin{align}
\Phi_X^{-1}&=\begin{bmatrix} \Phi_{oo}(z) & \Phi_{oh}(z) \\ \Phi_{ho}(z) & \Phi_{hh}(z) \end{bmatrix}^{-1}=\begin{bmatrix} J_{oo}(z) & J_{oh}(z) \\ J_{ho}(z) & J_{hh}(z) \end{bmatrix}\label{eqn:matrix_inversion}\\
&= (\mathbb{I}-\mathbb{H}(z))^{*}\Phi_E^{-1}(\mathbb{I}-\mathbb{H}(z))\nonumber.
\end{align}
\normalsize
Using the matrix inversion lemma \cite{horn2012matrix} in (\ref{eqn:matrix_inversion}), it follows that, $\Phi_{oo}^{-1} = J_{oo} - J_{oh}J_{hh}^{-1}J_{ho}=:\Gamma+\Delta+\Sigma,$ where,
\vspace{-0.3cm}
\small
\begin{align}\label{eqn:defsOfComps}
&\Gamma=(\mathbb{I}-\mathsf{H}_{oo}^{*})\Phi_{E_o}^{-1}(\mathbb{I}-\mathsf{H}_{oo}), \  \Delta =\mathsf{H}_{ho}^{*}\Phi_{E_h}^{-1}\mathsf{H}_{ho}, \\ &\Sigma = -\Psi^*\Lambda^{-1}\Psi, \text{ where, }\Lambda = \mathsf{H}_{oh}^{*}\Phi_{E_o}^{-1}\mathsf{H}_{oh} + \Phi_{E_h}^{-1}, \nonumber \\
&\Psi = \mathsf{H}_{oh}^{*}\Phi_{E_o}^{-1}(\mathbb{I} - \mathsf{H}_{oo}) + \Phi_{E_h}^{-1}\mathsf{H}_{ho}. \nonumber
\end{align}
\normalsize

% We will show that the edge set constructed from the support structure of $\Phi^{-1}_{oo}$ results in a graph with many more spurious edges as compared to the setting discussed in the previous section where all nodes are observed. 

{The following Lemma \ref{lem:whenCompsAreZero} relates the $ \mathcal{G_T}$ to the structure of $\Phi^{-1}_{oo}.$}
\begin{lem}\label{lem:whenCompsAreZero}
The following assertions hold
\begin{enumerate}
\item \label{item:assumpGammaToBeZero} For a $i,j \in \mathcal{V}_o,$ if (i) there does not exist a $k \in \mathcal{V}_o\setminus \{i,j \}$ such that $i-k-j$ in $\mathcal{G}_{T}$ and (ii) $i-j$ is not in $\mathcal{G}_{T}$, then ${\boldsymbol{B}_{i}^{\prime}}\Gamma{\boldsymbol{B}_{j}}=\boldsymbol{0}_{T\times T}.$ 
\item \label{item:assumpDeltaToBeZero} For a $i,j\in \mathcal{V}_o,$ there does not exist a $l \in \mathcal{V}_h$ such that $i-l-j$ in $ \mathcal{G_T}$, then ${\boldsymbol{B}_{i}^{\prime}}\Delta{\boldsymbol{B}_{j}}=\boldsymbol{0}_{T\times T}.$
\item \label{item:assumpLambdaDiagonal} For a $l_1,l_2 \in \mathcal{V}_h,$ there does not exist a $i \in \mathcal{V}_o$ such that $l_1-i-l_2$ in $ \mathcal{G_T}$, then $\Lambda$ is block diagonal and Hermitian.
\item \label{item:assumpForPsiToBeZero} Suppose in $\mathcal{G}_{T},$ for $j \in \mathcal{V}_o$ and $l \in \mathcal{V}_h$; (i) $j-l$ is not present and (ii) there is no path of the form $j-p-l$ with $p\in \mathcal{V}_o\setminus{j}$, then ${\boldsymbol{B}_{l}^{\prime}}\Psi{\boldsymbol{B}_{j}}=\boldsymbol{0}_{T \times T}.$ 

\item \label{item:assumpForSigmaToBeZero}Suppose $\Lambda$ is block diagonal and Hermitian and if in $\mathcal{G}_{T}$, for $i,j \in \mathcal{V}_o$ and $l \in \mathcal{V}_h$, there are no paths of the form $i-p-l$ or $i-l$ and $j-p^*-l$ or $j-l$ for any $p\in \mathcal{V}_o\setminus{i}$ and $p^*\in \mathcal{V}_o \setminus{j}$, then ${\boldsymbol{B}_{i}^{\prime}}\Sigma{\boldsymbol{B}_{j}}=\boldsymbol{0}_{T \times T}.$
\end{enumerate}
\end{lem}

\begin{pf}
The proof is similar to the proof of Lemma $3.1$ from \cite{talukdar2015reconstruction}, but extended to complex-valued cyclostationary processes. Its important to note the salient points (i) for $i,j \in \mathcal{V}_o, \ \mathsf{H}_{ij}:=\mathsf{H}_{oo}(i,j) $ is a $T \times T$ transfer matrix and $\mathsf{H}_{ii}= \boldsymbol{0}$ (ii) $\Phi^{-1}_{E},\ \Phi^{-1}_{E_o},\ \Phi^{-1}_{E_h}$ are block diagonal Hermitian matrices of size $nT \times nT, \ mT \times mT$ and $(n-m)T \times (n-m)T,$ respectively.
\begin{enumerate}
 \item It follows from (\ref{eqn:defsOfComps}), 
 \small 
 \begin{align}\label{Observed}
 {\boldsymbol{B}_{i}^{\prime}}\Gamma{\boldsymbol{B}_{j}}=&-{\boldsymbol{B}_{i}^{\prime}}\Phi_{E_o}^{-1}\mathsf{H}_{oo}{\boldsymbol{B}_{j}}-{\boldsymbol{B}_{i}^{\prime}}{\mathsf{H}^{*}_{oo}}\Phi_{E_o}^{-1}{\boldsymbol{B}_{j}}\nonumber \\
 &+{\boldsymbol{B}_{i}^{\prime}}{\mathsf{H}^{*}_{oo}}\Phi_{E_o}^{-1}\mathsf{H}_{oo}{\boldsymbol{B}_{j}},\nonumber \\
 {\boldsymbol{B}_{i}^{\prime}}\Gamma{\boldsymbol{B}_{j}}=-&\Phi_{E_i}^{-1}\mathsf{H}_{ij}-{\mathsf{H}^{*}_{ji}}\Phi_{E_j}^{-1}+\sum_{k=1}^{m}{\mathsf{H}^{*}_{ki}}\Phi_{E_k}^{-1}\mathsf{H}_{kj}.
 \end{align}
 \normalsize
 Here, $i,j \in \mathcal{V}_o.$ 
 If $(j,i) \notin \mathcal{E}_o,$ then $\mathsf{H}_{ij}= \boldsymbol{0}_{T \times T}.$ Similarly, if $(i,j) \notin \mathcal{E}_o,$ then $\mathsf{H}_{ji}= \boldsymbol{0}_{T \times T}.$ 
 In $\mathcal{G}_o,$ if there does not exist $k\in \mathcal{V}_o\setminus \{i.j \},$ such that $\{(k,i),(k,j)\} \in \mathcal{E}_o,$ then the third term is $\boldsymbol{0}_{T \times T}.$ 

 \item From (\ref{eqn:defsOfComps}), it follows that,
 \small
 \begin{align*}
 {\boldsymbol{B}_{i}^{\prime}}\Delta {\boldsymbol{B}_{j}}&=\sum_{l\in \mathcal{V}_h}[{{\mathsf{H}_{ho}(l,i)}}]^{*}\Phi_{E_h}^{-1}(l,l)\mathsf{H}_{ho}(l,j),\\ 
 {\boldsymbol{B}_{i}^{\prime}}\Delta{\boldsymbol{B}_{j}}&=\sum_{l\in \mathcal{V}_h}[{{\mathsf{H}_{li}}}]^{*}\Phi_{E_l}^{-1}\mathsf{H}_{lj}.
 \end{align*}
 \normalsize
 
 For a given $i,j \in \mathcal{V}_o,$ if there does not exist a $l \in \mathcal{V}_h$, such that $\{ (l,i),(l,j)\} \in \mathcal{E}_h.$ 
 
 \item Consider two distinct hidden nodes $l_1,l_2 \in \mathcal{V}_h.$ If there does not exist $i \in \mathcal{V}_o,$ such that $\{(i,l_1),(i,l_2)\} \in \mathcal{E}_h,$ then $\mathsf{H}_{oh}(i,l_1)= \mathsf{H}_{il_1}= \boldsymbol{0}_{T \times T}$ and $\mathsf{H}_{oh}(i,l_2)= \mathsf{H}_{il_2}= \boldsymbol{0}_{T \times T}.$ Thus, from (\ref{eqn:defsOfComps}), it follows that ${\boldsymbol{B}_{l_1}^{\prime}}\Lambda{\boldsymbol{B}_{l_2}}=\sum_{i\in \mathcal{V}_o} [\mathsf{H}_{oh}(i,l_1)]^{*}\Phi^{-1}_{E_i}\mathsf{H}_{oh} (i,l_2) $ is $\boldsymbol{0}_{T \times T}$. The diagonal block of $\Lambda$ is given by ${\boldsymbol{B}_{l_1}^{\prime}}\Lambda{\boldsymbol{B}_{l_1}}=\sum_{i\in \mathcal{V}_o} {\mathsf{H}_{oh}(i,l_1)]^{*}}\Phi^{-1}_{E_i}\mathsf{H}_{oh} (i,l_1)+ \Phi^{-1}_{E_h}(l_1,l_1),$ which is Hermitian of size $T \times T.$ Thus, $\Lambda$ is block diagonal Hermitian. 
 \item Suppose $(l,j)\notin \mathcal{E}_{T_h},$ then $\mathsf{H}_{oh}(j,l)= \mathsf{H}_{ho}(l,j)= \boldsymbol{0}_{T \times T}.$ If there does not exist a $p \in \mathcal{V}_o\setminus{j},$ such that $\{(p,j),(p,l) \}\in \mathcal{E},$ then $\mathsf{H}_{oh}(p,l)= \mathsf{H}_{ho}(p,j)= \boldsymbol{0}_{T \times T}.$ Thus, from (\ref{eqn:defsOfComps}), the $(l,j)^{th}$ block of $\Psi$ is given by 
 ${\boldsymbol{B}_{l}^{\prime}}\Psi{\boldsymbol{B}_{j}} = {[\mathsf{H}_{oh}(j,l)]^{*}}\Phi_{E_j}^{-1} -\sum_{p=1}^{m}{[\mathsf{H}_{oh}(p,l)]^{*}}\Phi_{E_p}^{-1}\mathsf{H}_{oo}(p,j)+\Phi_{E_h}^{-1}(l,l)\mathsf{H}_{ho}(l,j) = \boldsymbol{0}.$

 \item For a given $i,j \in \mathcal{V}_o,$ and for any $l \in \mathcal{V}_h$ the following holds (i) $ {\boldsymbol{B}_{l}^{\prime}}\Psi{\boldsymbol{B}_{i}}= \boldsymbol{0}$ if there does not exists $i \leftarrow l,i \rightarrow l$ and $i \rightarrow p \leftarrow l$ for any $p \in \mathcal{V}_o\setminus{i}$ in $\mathcal{E},$ 
 (ii) $ {\boldsymbol{B}_{l}^{\prime}}\Psi{\boldsymbol{B}_{j}} = \boldsymbol{0}$ if there does not exists $j \leftarrow l,j \rightarrow l$ and $j \rightarrow p \leftarrow l$ for any $p \in \mathcal{V}_o\setminus{j}$ in $\mathcal{E}.$ Thus, ${\boldsymbol{B}_{i}^{\prime}}\Sigma{\boldsymbol{B}_{j}}= \boldsymbol{0}_{T \times T}.$
\end{enumerate}
\end{pf}

Since, we are restricting our attention to bidirected networks, we make the following assumption.
\begin{assum}\label{ass:bidirected}
If $\mathsf{H}_{ji}(z)\neq 0,$ then $\mathsf{H}_{ij}(z)\neq 0.$
\end{assum}

\begin{assum}\label{ass:4hops}
The hidden nodes in $ \mathcal{G_T}$ are at least four or more hops away from each other.
\end{assum}

\begin{thm}\label{thm:4.1}
Consider a linear dynamical system with topology $ \mathcal{G_T}$ such that Assumptions \ref{ass:bidirected}, \ref{ass:timeperiod} and \ref{ass:4hops} hold. Then ${\boldsymbol{B}_{i}^{\prime}}\Phi^{-1}_{oo}{\boldsymbol{B}_{j}} \neq \boldsymbol{0}$ for all $\omega \in [0, 2\pi)$, implies that, $i$ and $j$ are within four hops of each other in the graph $ \mathcal{G_T}$.
\end{thm}

\begin{pf}
Given that ${\boldsymbol{B}_{i}^{\prime}}\Phi^{-1}_{oo}{\boldsymbol{B}_{j}} \neq \boldsymbol{0},$ then either (i) ${\boldsymbol{B}_{i}^{\prime}}\Gamma{\boldsymbol{B}_{j}} \neq \boldsymbol{0}$ or (ii) ${\boldsymbol{B}_{i}^{\prime}}\Delta{\boldsymbol{B}_{j}} \neq \boldsymbol{0}$ or (iii) ${\boldsymbol{B}_{i}^{\prime}}\Sigma{\boldsymbol{B}_{j}} \neq \boldsymbol{0}.$ The proof is by enumerating all the possible cases as follows.\\
(i) ${\boldsymbol{B}_{i}^{\prime}}\Gamma{\boldsymbol{B}_{j}} \neq \boldsymbol{0}$ implies that either $i\leftarrow j$ or $i \rightarrow j$ or $i \rightarrow p \leftarrow j$ exists in $\mathcal{E}_o,$ for some $p \in \mathcal{V}_o \setminus\{i.j \}.$ This is evident from $1)$ of Lemma \ref{lem:whenCompsAreZero}. \\
(ii) From $2)$ of Lemma \ref{lem:whenCompsAreZero}, ${\boldsymbol{B}_{i}^{\prime}}\Delta{\boldsymbol{B}_{j}} \neq \boldsymbol{0}$ implies there exists a $l \in \mathcal{V}_h, $ such that $i \rightarrow l \leftarrow j$ exists in $\mathcal{E}_h.$ \\
(iii) From $3),\ 4)$ and $5)$ of Lemma $\ref{lem:whenCompsAreZero},$ ${\boldsymbol{B}_{i}^{\prime}}\Sigma{\boldsymbol{B}_{j}} \neq \boldsymbol{0},$ implies that there exists a hidden node $l \in \mathcal{V}_h$ such that (a) either $i \leftarrow l$ or $i \rightarrow l$ or $i \rightarrow p \leftarrow l$ exists in $\mathcal{E}$ for some $p \in \mathcal{V}_o\setminus{i}$ and (b) either $j \leftarrow l$ or $j \rightarrow l$ or $j \rightarrow p \leftarrow l$ exists in $\mathcal{E}$ for some $p \in \mathcal{V}_o\setminus{j}.$

In (i), (ii) and (iii), the nodes $i$ and $j$ are within four hops of each other in $ \mathcal{G_T}.$ 
\end{pf}

\begin{rem}\label{rem:leaf_hidden_distance}
From  proof of Theorem \ref{thm:4.1}, it follows that if an observable node $i$ is  more than 2 hops away from any hidden node, then, for all $j \in \mathcal{V}_o\setminus{i},$ the following holds: (i) ${\boldsymbol{B}_{i}^{\prime}}\Delta{\boldsymbol{B}_{j}}\ =\ {\boldsymbol{B}_{i}^{\prime}}\Sigma{\boldsymbol{B}_{j}}=\ \boldsymbol{0},$ and (ii) ${\boldsymbol{B}_{i}^{\prime}}\Phi^{-1}_{oo}{\boldsymbol{B}_{j}} \neq \boldsymbol{0}$ implies that ${\boldsymbol{B}_{i}^{\prime}}\Gamma{\boldsymbol{B}_{j}} \neq \boldsymbol{0}$ {(no hidden node contribution)}. 
\end{rem}

% \begin{comment}
% \begin{pf}
% It follows from (\ref{invPhi}) that, $\Phi_{oo}^{-1}(i,j)(\omega) = \Gamma(i,j)+\Delta(i,j)+\Sigma(i,j)$. %\\&= \mathcal{N}_{o,1} + \mathcal{N}_{o,2} + \mathcal{N}_{h,2} + \mathcal{N}_{h,3} + \mathcal{N}_{h,4}.
% %\end{align*}
% Suppose $i$ and $j$ are more than four hops away. We will conclude that $\Phi_{oo}^{-1}(i,j)(\omega) = 0$ almost surely.

% As two and one hop paths are not present, it follows from \ref{item:assumpGammaToBeZero}) and \ref{item:assumpDeltaToBeZero}) of Lemma~\ref{lem:whenCompsAreZero} that $\Gamma(i,j)=\boldsymbol{0}$ and $\Delta(i,j)=\boldsymbol{0}$ respectively.

% From the assumption that unobserved nodes are at least four or more hops away, it follows from \ref{item:assumpLambdaDiagonal}) of Lemma~\ref{lem:whenCompsAreZero} that $\Lambda$ is block diagonal and Hermitian. 

% Suppose there are paths of the form $i-p-k$ and $j-p^*-k$ then $i-p-k-p^*-j$ is a four hop path that connects $i$ and $j$, which contradicts that $i$ and $j$ are more than four hops away. Thus, paths of the form $i-p-k$ and $j-p^*-k$ cannot be present simultaneously. Similarly, one can show that paths of the form $i-k$ and $j-p^{*}-k$ or $i-p-k$ and $j-k$ cannot be present as it would imply $i,j$ are three hop neighbors. Thus, from \ref{item:assumpForSigmaToBeZero}) of Lemma~\ref{lem:whenCompsAreZero}, we conclude that $\Sigma(i,j)=\boldsymbol{0}.$ This implies that $\Phi_{oo}^{-1}(i,j)=\boldsymbol{0}$ and completes the proof.
% \end{pf}
% \end{comment}

An instance of the transfer function matrix $\mathbb{H}(z)$ can be designed so that converse of Theorem \ref{thm:4.1} does not hold. However, such cases are pathological and the converse of Theorem \ref{thm:4.1} holds for all practical purposes. Thus, based on Theorem \ref{thm:4.1}, we construct an undirected graph based on the structure of $\Phi^{-1}_{oo}$ as follows: Initialize $\mathcal{E}_c$ as $\{ \}$ and set $\mathcal{V}_o$ as $\{1,2,\cdots,m \}.$ Construct the undirected edge set {$\mathcal{E}_c:= \{(i,j)|\ {\boldsymbol{B}_{i}^{\prime}}\Phi^{-1}_{oo}{\boldsymbol{B}_{j}} \neq \boldsymbol{0} \}.$} The undirected graph $\mathcal{G}_c:=(\mathcal{V}_o, \mathcal{E}_c)$ constitutes edges between the nodes that are within four hops of each other in $ \mathcal{G_T}.$
From Theorem \ref{thm:4.1}, it follows that $\mathcal{E}_{T_o} \subset \mathcal{E}_c,$ and hence all the true edges among the observed nodes can be recovered in addition to spurious edges. We use $\mathcal{G}_c(\mathcal{V}_o, \mathcal{E}_c)$ in the following theorems to obtain $\mathcal{G}_{T_o}.$ We identify $\mathcal{E}_{T_o}$, number of hidden nodes, and their neighbors solely from $\mathcal{G}_c(\mathcal{V}_o, \mathcal{E}_c).$ We will make the following assumption of $ \mathcal{G_T}$ for recovering $\mathcal{E}_{T_o}$, number of hidden nodes, and their neighbors from $\mathcal{G}_c.$

\begin{assum}\label{ass:tree}
$ \mathcal{G_T}$ is {tree}, that is for any two nodes $i,j \in \mathcal{V},$ there is a unique path connecting nodes $i$ and $j$ in $ \mathcal{G_T}.$ Further, every hidden node is at least three hops away from all leaf nodes in $ \mathcal{G_T}$.
\end{assum}

Note that with Assumption \ref{ass:tree}, $ \mathcal{G_T}$ does not posses cycles and henc, the Assumption \ref{ass:atmost1} holds. Thus, Theorem \ref{thm:3.2} and Theorem \ref{thm:3.3} from Appendix \ref{app:thm3.3_proof} holds. The nodes with degree $1$ are termed as leaf nodes $V_l$, while other nodes in $\mathcal{V}$ are referred to as non-leaf nodes $V_{nl}$. By Assumption \ref{ass:tree}, $V_l \in \mathcal{V}_o.$

\subsection{Reconstructing $\mathcal{G}_{T_o}$ from $\mathcal{G}_c$}
Given that $ \mathcal{G_T}$ satisfies Assumption \ref{ass:4hops} and Assumption \ref{ass:tree}, we propose an algorithm for identifying the leaf and non-leaf observable nodes, and the topology restricted to observable nodes. The phase result of the eigenvalues developed in the previous section is not applicable to identify the topology restricted to observed nodes as we do not know the locations of the hidden nodes. Here we exploit graphical separation and the {tree }topology to identify the true edges from $\mathcal{G}_c$. Consider an undirected graph $U(V_u,E_u),$ and nodes $i,j \in V_u.$ The vertex set $Z\in V_u \setminus\{i,j\}$ is said to separate nodes $i$ and $j$ in $U$ if all the paths from $i$ to $j$ in $U$ contains at least a node in $Z.$ If $Z$ separates $i$ and $j$ in $U,$ then we say $sep(i,Z,j).$ The following theorem enables us to identify the non-leaf nodes and edges among them from $\mathcal{G}_c$. 

\begin{thm}\label{thm:4.2}
Consider a networked LDS (\ref{eqn:latent_LDM}), such that Assumptions \ref{ass:bidirected}, \ref{ass:timeperiod}, \ref{ass:4hops} and \ref{ass:tree} hold. There exist distinct nodes $a,b,c,d \in \mathcal{V}_o$ such that $sep(c,\{a,b\},d)$ holds in $\mathcal{G}_c$ if and only if $(a,b) \in  \mathcal{E_T}$ and $a,b$ are non-leaf nodes.
\end{thm}

\begin{pf}
Suppose $a-b$ is not an edge in $ \mathcal{G_T}$. Let $p:= c -\pi_{h,1}- \pi_1 - \pi_2 - \pi_{h,2}-\pi_3- \cdots-\pi_m-\pi_{h,j}- d$ be the unique path between $c$ and $d$ in $ \mathcal{G_T}$ such that $\{\pi_1,\pi_2,...,\pi_m\}$ are observed nodes and $\{\pi_{h,1},\pi_{h,2},...,\pi_{h,j}\}$ are hidden nodes. There are three possibilities - (i) neither of $a,b$ belong to $\{\pi_1,...,\pi_m\}$, (ii) either $a$ or $b$ but not both belong to $\{\pi_1,...,\pi_m\}$ and (iii) both $a,b$ belong to $\{\pi_1,...,\pi_m\}$ with $a-b$ not being
an edge.

\noindent (i) If $a$ and $b$ do not belong to $\{\pi_1,\cdots,\pi_m\}$.
Then $c - \pi_1 - \pi_2 - \cdots - \pi_m - d$ is a path in $\mathcal{G}_c$ with no intermediate node in the path being $a$ or $b$. Thus, $sep(c,\{a,b\},d)$ does not hold, which is a contradiction.

\noindent (ii) If $a$ belongs to $\{\pi_1,\cdots,\pi_m\}$ but not $b$. Let $\pi_k = a$.
Then $c - \pi_1 - \pi_2 - \cdots - \pi_{j} - a - \pi_{l} \cdots - \pi_m - d$ is a path in $\mathcal{G}_c$ which is not separated by $\{a,b\}$. Thus, $sep(c,\{a,b\},d)$ does not hold which is a contradiction. Similarly, one can arrive at a contradiction for the case $b$ belongs to $\{\pi_1,\cdots,\pi_m\}$ but not $a$.

\noindent (iii) If both $a$ and $b$ belong to $\{\pi_1,\cdots,\pi_m\}$, and $a-b$ is not in $ \mathcal{G_T}$. Let $a = \pi_e$ and $b = \pi_j$ Then

(a) $a,b$ are two-hop neighbors through an observed node $\pi_g$.

Let $\pi_d,\pi_l$ be an observed neighbor of $a,b,$ respectively and $\pi_g$ be the common neighbor of $a,b$. Then,
$c-\pi_1-...-\pi_d-\pi_g-\pi_l-...-\pi_m-d$ is a path in $\mathcal{G}_c$ which is not separated by $\{a,b\}$. Thus, $sep(c,\{a,b\},d)$ does not hold and is a contradiction.

(b) $a,b$ are two-hop neighbors through an unobserved node $\pi_{h,g}$.

Let $\pi_d,\pi_l$ be an observed neighbor of $a,b,$ respectively and $\pi_{h,g}$ is the common unobserved neighbor of $a,b$. Then,
$c-\pi_1-...-\pi_d-\pi_l-...-\pi_m-d$ is a path in $\mathcal{G}_c$ which is not separated by $\{a,b\}$. Thus, $sep(c,\{a,b\},d)$ does not hold and is a contradiction.

(c) Let $a,b$ are three hop neighbors with one hidden node($\pi_{h,f}$) and one observed node($\pi_{f}$) in between $a,b$.

Let $\pi_d$ be neighbor of $a$ on the other side of $b$. Similarly, let $\pi_l$ be another neighbor of $b$ in the other direction of $a$. Then, $c-\pi_1-...-\pi_d-a-\pi_f-\pi_{h,f}-b-\pi_l-...-\pi_m-d$ is a path in $ \mathcal{G_T}$, with $c-\pi_1-...-\pi_d-\pi_f-\pi_l-...-\pi_m-d$ being a path in $\mathcal{G}_c$ not separated by $\{a,b\}$ and is a contradiction.

(d) Let $a,b$ are four hop neighbors with $\pi_f, \pi_g$ being observed neighbors of $a$ and $b,$ respectively and $\pi_{h,f}$ being an unobserved neighbor of $\pi_f,\pi_g$ in $ \mathcal{G_T}$.

The path in $ \mathcal{G_T}$ is of the form $c-\pi_1-...-\pi_d-a-\pi_f-\pi_{h,f}-\pi_{g}-b-\pi_l-...-\pi_m-d$. Then $c-\pi_1-...-\pi_d-\pi_f-\pi_g-\pi_l-...-\pi_m-d$ is a path in $\mathcal{G}_c$ which is not separated by $a,b$ and is a contradiction.

(e) Let $a,b$ be more than four hops away such that $a-\pi_f-...-\pi_{h,f}-...-\pi_g-...\pi_{h,g}-...-\pi_h-...-\pi_l-b$. Using the same reasoning as before one can show that a path exists in $\mathcal{G}_c$ which does not contain both $a$ and $b$, that is, there exist a path which is not separated by $\{a,b\}$ in $\mathcal{G}_c$. Thus, $sep(c,\{a,b\},d)$ does not hold in $\mathcal{G}_c$ and is a contradiction.

\noindent As all cases have been exhausted we conclude that $sep(c,\{a,b\},d)$ in $\mathcal{G}_c$ is not possible, which is a contradiction. Hence, $a-b$ is a true edge in $ \mathcal{G_T}$. Both $a,b$ have degree at least two as they have at least two neighbors, hence, are non-leaf nodes. This proves the theorem.
\end{pf}

The conclusion of Theorem \ref{thm:4.2} is that, if $(a,b)$ is a spurious edge between observable non-leaf nodes $a,b \in \mathcal{V}_o$, then there exist no $c,d$ different from $a,b$ such that $sep(c,\{a,b\},d)$ holds in $\mathcal{G}_c$. This provides a graph based test to identify the true edges among non-leaf observable nodes from $\mathcal{G}_c.$ Thus, all the non-leaf nodes $V_{nl} \in \mathcal{V}_o,$ the edges $\{(i,j)| i,j \in V_{nl}, (i,j) \in  \mathcal{G_T} \}$ are recovered. Since a leaf node is not hidden by Assumption \ref{ass:tree}, the set of leaf nodes are given by $V_l:= \mathcal{V}_o \setminus V_{nl}.$ Suppose $l \in V_l,$ then it has a single non-leaf neighbor in $ \mathcal{G_T},$ since degree of $l$ is one. Based on Assumption \ref{ass:tree}, the node $l$ is at least three hops away from any hidden node in $ \mathcal{G_T}$. From Lemma \ref{lem:whenCompsAreZero} and Remark \ref{rem:leaf_hidden_distance}, the spurious edges connected with $l$ in $\mathcal{G}_c$ include those up to its two-hop neighbors in $ \mathcal{G_T}$. {Thus, the two-hop neighbors of a leaf node in $ \mathcal{G_T}$ should be identified and pruned. The following theorem is helpful for identification of two-hop neighbors of a leaf node.}

\begin{thm}\label{thm:4.3}
Consider a networked LDS described by (\ref{eqn:latent_LDM}) such that Assumptions \ref{ass:transferFunction_property}, \ref{ass:timeperiod}, \ref{ass:bidirected}, \ref{ass:4hops} and \ref{ass:tree} hold. Let $i \in V_l,$ and $j \in V_{nl}$ be a neighbor of $i$ in $\mathcal{G}_c$. Suppose $\Eps_{ij} \in \C^{T \times 1}$ be the eigenvalues of ${\boldsymbol{B}_{i}^{\prime}}\Phi^{-1}_{oo}{\boldsymbol{B}_{j}}.$ Then, $\phase{\Eps_{ij}(l)}$ is a constant for $l\in \{1,2,\cdots,T \}$ and all $\omega \in [0, 2 \pi)$ if and only if $i,j$ are two-hop neighbors in $ \mathcal{G_T}$.
\end{thm}

\noindent{\bf PROOF:}
The proof follows directly from Theorem \ref{thm:3.3}.

Based on Theorem \ref{thm:4.2} and \ref{thm:4.3}, we propose an Algorithm $2$ which reconstructs $\mathcal{G}_{T_o}.$ The Algorithm $2$ is shown in Appendix \ref{app:alg2}. As asserted earlier, the two-hop neighbors of a leaf node are identified using Theorem \ref{thm:4.3} when Assumption \ref{ass:transferFunction_property} holds. When Assumption \ref{ass:transferFunction_property} does not hold, the spurious edge associated with a leaf nodes cannot be identified, and thus, $\mathcal{E}_{T_o} \subset \mathcal{E}_{\overline{\mathcal{T}}}.$ $\mathcal{E}_{\overline{\mathcal{T}}} \setminus \mathcal{E}_{T_o}$ is the set of two-hop neighbors of the leaf nodes. After learning $\mathcal{G}_{T_o},$ we find $\mathcal{V}_h$ and $\mathcal{E}_{T_h}$ using the next theorem (note that the theorem statement appeared in a preliminary conference paper \cite{talukdar2018topology} without proof).

\begin{thm}\label{thm:insert_hidden}
Consider a networked LDS described by (\ref{eqn:latent_LDM}) such that Assumptions \ref{ass:timeperiod}, \ref{ass:bidirected}, \ref{ass:4hops} and \ref{ass:tree} hold. Suppose $\overline{\mathcal{T}}_1,\overline{\mathcal{T}}_2$ are two disconnected components in $\mathcal{G}_{T_o}(\mathcal{V}_o,\mathcal{E}_{T_o})$ with observed nodes $c \in \overline{\mathcal{T}}_1$ and $e \in \overline{\mathcal{T}}_2$. If for all $ b\in \overline{\mathcal{T}}_1$, for all $f\in \overline{\mathcal{T}}_2$ where $b-c$ and $e-f$ are edges in true topology $\mathcal{G}_{T_o}$ and $b,c,e,f$ form a clique in $\mathcal{G}_c$, then there exists a $d \in \mathcal{V}_h$ such that $c-d-e$ is a path in $ \mathcal{G_T}.$ 
\end{thm}

\begin{pf}
Since, $ \mathcal{G_T}$ satisfies Assumptions \ref{ass:4hops}, it follows that $\mathcal{G}_{T_o}(\mathcal{V}_o,\mathcal{E}_{T_o})$ is a union of disconnected connected components, where each component has at least three nodes. Suppose $\overline{\mathcal{T}}_1,\overline{\mathcal{T}}_2$ are two disconnected components in $\mathcal{G}_{T_o}(\mathcal{V}_o,\mathcal{E}_{T_o})$ with observed nodes $c \in \overline{\mathcal{T}}_1$ and $e \in \overline{\mathcal{T}}_2$. Consider $b \in \overline{\mathcal{T}}_1, b \neq c$ and $f \in \overline{\mathcal{T}}_2, f \neq e,$ such that $b-c$ and $e-f$ exists in $ \mathcal{G_T}.$ 
From Assumptions \ref{ass:4hops} and \ref{ass:tree}, there exists an observable node $a \in \overline{\mathcal{T}}_1, g \in \overline{\mathcal{T}}_2$ such that $a - b - c$ exists and $e-f-g$ exists in $ \mathcal{G_T}.$ 

We will show that if there is no $l \in \mathcal{V}_h$ such that $c-l-e$ exists in $ \mathcal{G_T},$ then $b,c,e,f$ cannot form a clique in $\mathcal{G}_c.$ If there is no such $l,$ then $a-b-c-d_1- d - d_2 - e-f-g$ exists in $ \mathcal{G_T}$ for some $d_1, d_2 \in \mathcal{V}_o$ and $d \in \mathcal{V}_h.$ Node $d$ exists because $a,b,c \in \overline{\mathcal{T}}_1$
and $e,f,g \in \overline{\mathcal{T}}_2.$ From Lemma \ref{lem:whenCompsAreZero}, it is evident that $b,c,d,e$ does not form a clique in $\mathcal{G}_c,$ which is a contradiction. Same holds if one $d_1$ or $d_2$ is present. Thus there does not exists $d_1$ and $d_2,$ such that $a-b-c-d_1- d - d_2 - e-f-g$ exists in $ \mathcal{G_T}.$ Therefore, $b-c-d - e-f$ exists in $ \mathcal{G_T}.$
\end{pf}

% \small
% \begin{algorithm}
% \small{
% \caption{Reconstructing $\mathcal{V}_h$ and $\mathcal{E}_{T_h}$}
% \textbf{Input:} $\overline{\mathcal{T}}=(\mathcal{V}_o,\mathcal{E}_{\overline{\mathcal{T}}})=\cup_{j=1}^{h}\overline{\mathcal{T}}_j$\\
% \textbf{Output:} $\tilde{\mathcal{T}} = (V_{\tilde{\mathcal{T}}},\mathcal{E}_{\tilde{\mathcal{T}}})$.
% \begin{algorithmic}[1]
% \State Node set $V_{\tilde{\mathcal{T}}} \gets \mathcal{V}_o$, edge set $\mathcal{E}_{\tilde{\mathcal{T}}}\gets \mathcal{E}_{\overline{\mathcal{T}}}$
% \State $h\gets$ Number of disjoint subgraphs in $\mathcal{G}_{T_o}$
% \ForAll{$j \in \{1,2,...,h\}$}
% \ForAll{$i \in \{j+1,...,h\}$}
% \If{there exist a pair of nodes $a,b$ such that $a \in \overline{\mathcal{T}}_j$ and $b \in \overline{\mathcal{T}}_i$ such that all their neighbors in $\overline{\mathcal{T}}$ are connected in $\mathcal{G}_c$}
% \State $V_{\tilde{\mathcal{T}}} \gets V_{\tilde{\mathcal{T}}} \cup l_j$
% \State $\mathcal{E}_{\tilde{\mathcal{T}}} \gets \mathcal{E}_{\tilde{\mathcal{T}}} \cup \{(a,l_j),(l_j,b)\}$
% \EndIf
% \EndFor\label{step1_b}
% \EndFor
% \State Merge hidden nodes that are neighbors of the same observed node.
% \end{algorithmic}}
% \end{algorithm}
% \normalsize

Let $\mathcal{G}_{T_o}$, the graph topology of observed nodes, be reconstructed using Algorithm $2$. As alluded earlier, $\mathcal{G}_{T_o}$ is a union of disconnected subgraphs. Let $h$ be the number of disconnected subgraphs in $\mathcal{G}_{T_o}.$ In Algorithm $3$, shown in Appendix \ref{app:alg3}, for each pair of disconnected subgraphs $\overline{\mathcal{T}}_i, \overline{\mathcal{T}}_j$, Theorem \ref{thm:insert_hidden} is checked to identify if a hidden node exists between them (Step $5$). If yes, then a hidden node is inserted and the edges to its neighbors are added. This completes the topology reconstruction of $\tilde{\mathcal{T}} = (V_{\tilde{\mathcal{T}}},\mathcal{E}_{\tilde{\mathcal{T}}}),$ which is identical to $ \mathcal{G_T}(\mathcal{V}_o \cup \mathcal{V}_h,  \mathcal{E_T}).$ 
Note that Theorem \ref{thm:insert_hidden} does not require Assumption \ref{ass:transferFunction_property} to hold. Consider $\tilde{\mathcal{T}}$ when Assumption \ref{ass:transferFunction_property} does not hold. The spurious edges associated with the leaf nodes are recovered in $\mathcal{E}_{\tilde{\mathcal{T}}},$ along with the true edges of the leaf nodes. The hidden nodes and their one-hop, two-hop neighborhoods in $\mathcal{G_T}$ are identified exactly. Moreover, if a hidden node is present at a $3-$ hop distance from a leaf node in $\mathcal{G_T}$, then the spurious edge associated with that leaf node in $\tilde{\mathcal{T}}$ 
can be identified from the recovered exact two-hop neighborhood of the hidden node. 
For a network with directed loops or loopy topology, topology learning with latent nodes for cyclostationary processes is a challenging task and future effort. In the next section, we provide numerical results to validate the theoretical algorithms presented in Section \ref{sec:three} and Section \ref{sec:four}.
\section{Results} \label{sec:five}
The illustrations presented in this section are based on networks with dynamic links among agents. CVX \cite{grant2014cvx} is used to solve the vector optimization (\ref{eqn:vec_opt_sparse}) and (\ref{eqn:RWFE}). 

%\subsection{Ocean Current Reconstruction from SST data}
\noindent \underline{\it Ocean Current Reconstruction from SST data: }
Ocean currents are bulk movement of water from one region to another in the ocean and has a network representation. The regions are considered as nodes of the network and the current direction is represented as a network edge. The currents are reconstructed using monthly mean COBE-SST2 Sea Surface Temperature alone, downloaded from \cite{cobe2sst_dataset}.  
%\cite{hirahara2014centennial}.
The nodes considered for the network reconstruction are spatially distant in the Atlantic region to begin with, and from the SST data of these nodes the edges of the network are reconstructed. The SST monthly mean data for the period 1850 -2019 consists of $2040$ data samples. {\color{red}}  $X(f)$  the FFT of the SST data at frequencies $f =\frac{2\pi k}{2040}, k = \{1,\cdots, \frac{2040}{2}-1 \}$ shows  a large magnitude  at 12 months, second highest at  6 months; there are  components  corresponding to 3 months and lower exists, but their magnitude is small.} 
%Fig. \ref{fig:SST_timeperiod} shows the power spectrum, which consists of peaks at the highlighted frequencies.
The time-period of the SST data is estimated to be 12 months, hence $T = 12.$ The values of $n,N, \lambda$ in solving (\ref{eqn:RWFE}) are $42 ,4, $ and $1.7,$ respectively. The threshold used in Algorithm $1$ is $0.04.$
Fig. \ref{fig:SST_reconstruction} shows the reconstructed currents along the west and east coast of Atlantic. The reconstructed currents in Fig. \ref{fig:SST_reconstruction} can be identified with known ocean currents:
gulf stream: $15 \rightarrow 14 \rightarrow 12 \rightarrow 10 \rightarrow 8 \rightarrow 17,$  Benguela current: $7 \rightarrow 4 ,$ Braziliean current: $1 \rightarrow 3,$ Malvinas current: $6 \rightarrow 5,$ Canary current: $9 \rightarrow 11 \rightarrow 13 \rightarrow 16.$  We validated our Algorithm against a known structure of the currents in North and South Atlantic. {\it It is important to note that the data generative model is unknown to us and more importantly, it is a non-linear model which does not conform to the generative model (linear) presented in the article. However, the algorithm can be applied on the data, and we showcase its performance by comparing with the true current direction from the literature.} The sample complexity analysis on the regularized Wiener filter estimator (\ref{eqn:RWFE}) is unknown. When the cyclostationary process time-period is one, $T = 1$ (wide-sense stationary), then the sample-complexity of the optimization problem (\ref{eqn:RWFE}) is solved in \cite{doddi2022efficient}. For theoretical guarantees, the choice of the regularization parameter, threshold, sample size $n, N$ need to be determined for $T>1$, and forms a part of our future work.

\begin{figure}[tb]
	\centering
	\begin{tabular}{c}
\includegraphics[width=0.75\columnwidth]{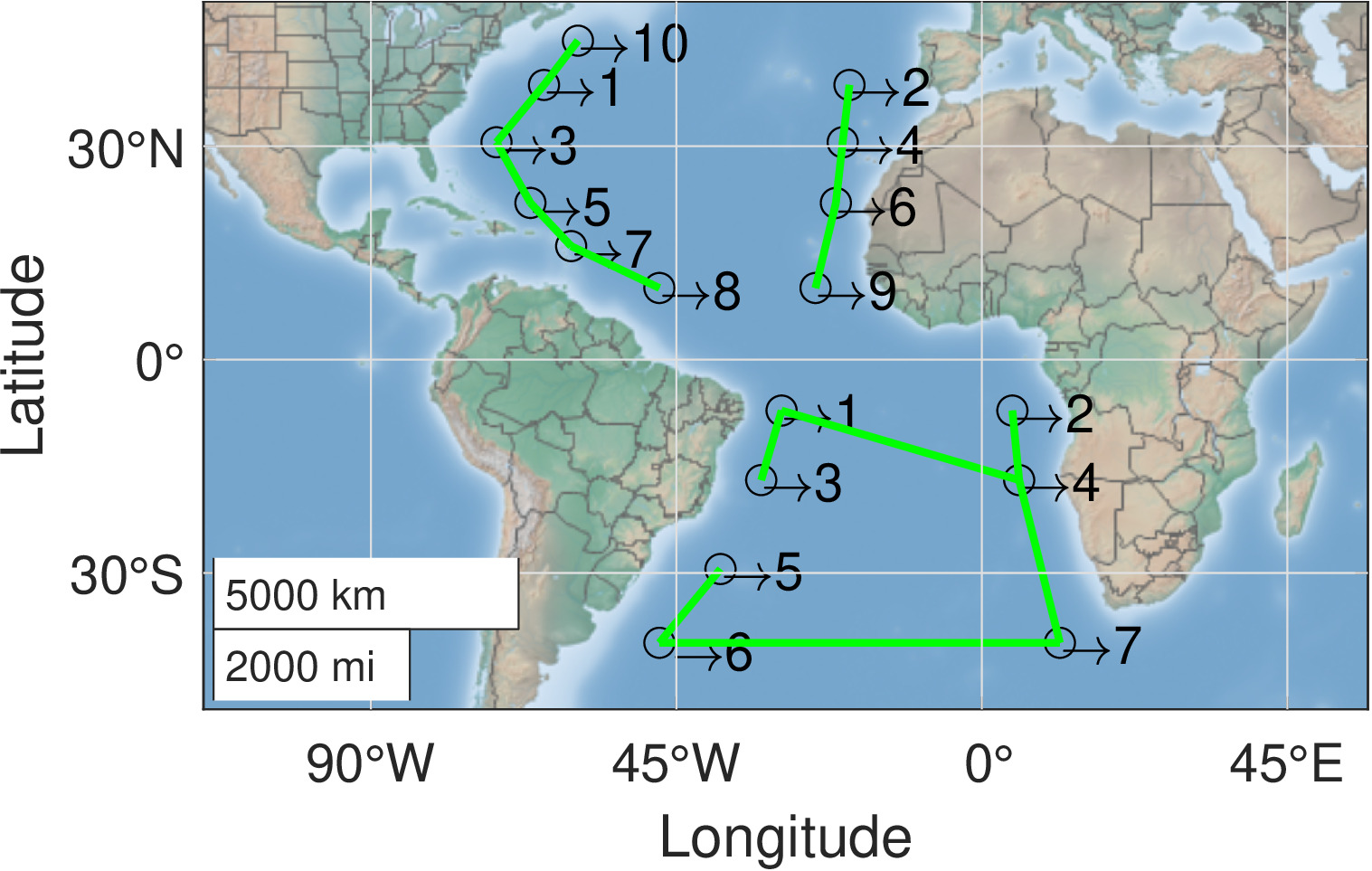} 
	\end{tabular} 
	\caption{Current reconstruction using 
 % the regularized Wiener filter estimated based on 
 frequency-domain regression presented in Subsection~\ref{subsec:RWFE}. Reconstructed currents are:
	gulf stream: $15 \rightarrow 14 \rightarrow 12 \rightarrow 10 \rightarrow 8 \rightarrow 17,$  Benguela current: $7 \rightarrow 4 ,$ Braziliean current: $1 \rightarrow 3,$ Malvinas current: $6 \rightarrow 5,$ Canary current: $9 \rightarrow 11 \rightarrow 13 \rightarrow 16.$ \label{fig:SST_reconstruction}}
\end{figure}

\subsection{Reconstruction of a Test network} 
We presented consistent topology reconstruction under full and partial observability for a test network containing $50$ nodes with $5$ latent nodes.

%\subsection{Test Network: Reconstruction under full observability:}
\noindent \underline{\it Test Network:  Reconstruction under full observability}
{We consider a generative graph consisting of observable $50$ nodes, with topology $ \mathcal{G_T}(\mathcal{V},  \mathcal{E_T})$ shown in Fig. \ref{fig:Results_topology}(a). $ \mathcal{E_T}$ includes green colored edges. All the edges (including green colored) are bidirected in the generative graph and hence posses directed loops,
but satisfy Assumption \ref{ass:atmost1}.
Note that generative graph is not restricted to bidirected network; the generative graph is allowed to be any directed graph satisfying Assumption \ref{ass:atmost1}. The time-series $x_i(k)$ is available at each node $i \in \{1,2,\cdots,50 \}.$ The dynamic links $h_{ji}$ in (\ref{eqn:cyclo_LDM}) are chosen to be FIR transfer functions with {\it complex-valued filter coefficients} satisfying Assumption \ref{ass:transferFunction_property}.} The exogenous input $e_1(k)$ is generated in MATLAB as $e_1(k) = cos(\pi k)w_1(k),$ where $w_1(k)$ is a zero mean wide sense stationary process. $e_i(k)$ is a cyclostationary process of period $T_1=2.$ The remaining exogenous inputs are mutually uncorrelated wide sense stationary processes, uncorrelated with $e_1(k).$ Using $\{e_i(k) \}_{i=1}^{50},$ time-series $\{x_i(k) \}_{i=1}^{50}$ are generated using (\ref{eqn:cyclo_LDM}). Each nodal time-series contains $628400$ samples. The simulations are done using MATLAB. 
% Although the theoretical guarantees for exact topology learning provided by all the algorithms discussed previously are in asymptotic sample limit, here we demonstrate that the error in topology reconstruction is significantly smaller for finite sample size and reduces when more data samples are available. 
Based on the periodogram analysis of the $\{x_i(k)\}_{i=1}^{50},$ we get  $T=2$. Each $x_i(k)$ is lifted to a vector time-series $X_i(k)=[x_i(2k), x_i(2k+1)]^{\prime}$ for $i=\{1,\cdots,50\}$. Solve (\ref{eqn:vec_opt_sparse}) for all $j = \{1,\cdots,50 \}$ with a predetermined constants $\gamma = 0.1$ and $L_* =3.$ The moral graph ($\bar{\mathcal{E}}_M$) is reconstructed based on the condition  $H_{\infty}[ \boldsymbol{\hat{\mathsf{W}}}_{ji} ]+H_{\infty}[ \boldsymbol{\hat{\mathsf{W}}}_{ij}] > \tau,$ for $i,j \in \{ 1,\cdots,50 \}$ and $\tau=0.03.$ Fig. \ref{fig:Results_topology}(b) shows a color map of $H_{\infty}[ \boldsymbol{\hat{\mathsf{W}}}_{ji} ]+H_{\infty}[ \boldsymbol{\hat{\mathsf{W}}}_{ij}].$ Some of the spurious edges have very low color intensity and are not identified in \ref{fig:Results_topology}(b), which reduces the computational effort of pruning them. The values of threshold $\tau,\ \gamma$ and $L_*$ are tuned to recover the exact topology when the number of samples per node is $6(10)^{5}$. An edge $(j,i)  \in \bar{\mathcal{E}}_M$ is considered as spurious (strict spouses) if $0.5[|\phase{\Eps_{ji} (\omega)}|+|\phase{\Eps_{ij} (\omega)}|]$ is constant. Here, $\Eps_{ji}(\omega)$ is the vector of eigenvalues of $\boldsymbol{\hat{\mathsf{W}}}_{ji}(\omega).$ Such edges are pruned. $ \mathcal{G_T}$ is exactly recovered after pruning the spurious edges from $\bar{\mathcal{E}}_M.$

%\subsection{Topology reconstruction under partial observability}
\noindent\underline{\it Reconstruction under partial observability:}
{
Consider the generative graph consisting of $50$ nodes with topology $ \mathcal{G_T}$ shown in Fig. \ref{fig:Results_topology}(a) without green colored edges. The red colored nodes $\{46,47,48,49,50 \}$ are latent. We aim to reconstruct $ \mathcal{G_T}$ using the time-series $\{x_i(k) \}_{i=1}^{{45}}.$
The efficacy of the algorithms was assessed using comparisons with the truth, determined by the exact expression of power spectral density using the generative model. $T=2$ based on the periodogram analysis of $\{x_i(k)\}_{i=1}^{{45}}.$ The lifted processes are $X_i(k)=[x_i(2k), x_i(2k+1)]^{\prime}$ for $i=\{1,\cdots,{45}\}$. 
We solve (\ref{eqn:vec_opt_sparse}) for all $j = \{1,\cdots,45 \}$ with a predetermined constants $\gamma = 0.1$ and $L_* =3.$ Construct an edge set $\mathcal{E}_c$ by placing an undirected edge between two distinct nodes $j$ and $i$ if $H_{\infty}[ \boldsymbol{\hat{\mathsf{W}}}_{ji} ]+H_{\infty}[ \boldsymbol{\hat{\mathsf{W}}}_{ij}] > \tau,$ for $i,j \in \{ 1,\cdots,50 \}$ and  $\tau=0.03.$ The values of $\tau,\ \gamma$ and $L_*$ are tuned based on large sample size. Fig. \ref{fig:Results_Hnorm_Partial} shows a color map of $H_{\infty}[ \boldsymbol{\hat{\mathsf{W}}}_{ji} ]+H_{\infty}[ \boldsymbol{\hat{\mathsf{W}}}_{ij}].$ Note that $\mathcal{G}_c$ is not a moral graph of $\mathcal{G}$. We now apply Algorithm $2$ on $\mathcal{G}_c$ to obtain   $\mathcal{G}_{T_o}(\mathcal{V}_o,\mathcal{E}_{\bar{T}}).$ $\mathcal{G}_{T_o}$ matched with $\mathcal{G}_{T_o}.$ For illustrative purpose, we show phase of eigenvalues (see Fig. \ref{fig:Results_phase}) for various edges after applying the graphical separation in Algorithm $2.$ It is evident that edges $\{(14,45),(44,15),(30,32),(27,29)\}$ are spurious because their eigenvalue phase is constant for all frequencies. Similarly other spurious edges are eliminated to obtain $\bar{\mathcal{T}}.$ The reconstructed graph $\bar{\mathcal{T}}$ recovers $\mathcal{G}_{T_o}.$ Using Algorithm $3,$ the presence of hidden nodes and their neighbors are reconstructed. The reconstructed topology recovers the true topology shown in Fig. \ref{fig:Results_topology} (a).

\begin{figure}[hb]
	\centering
	\begin{tabular}{cc}
	\includegraphics[width=0.4\columnwidth]{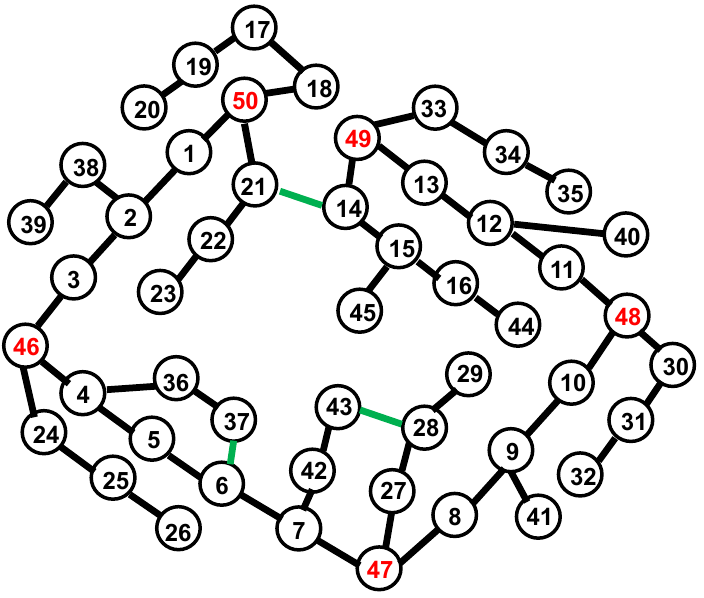} &
 	\includegraphics[width=0.5\columnwidth]{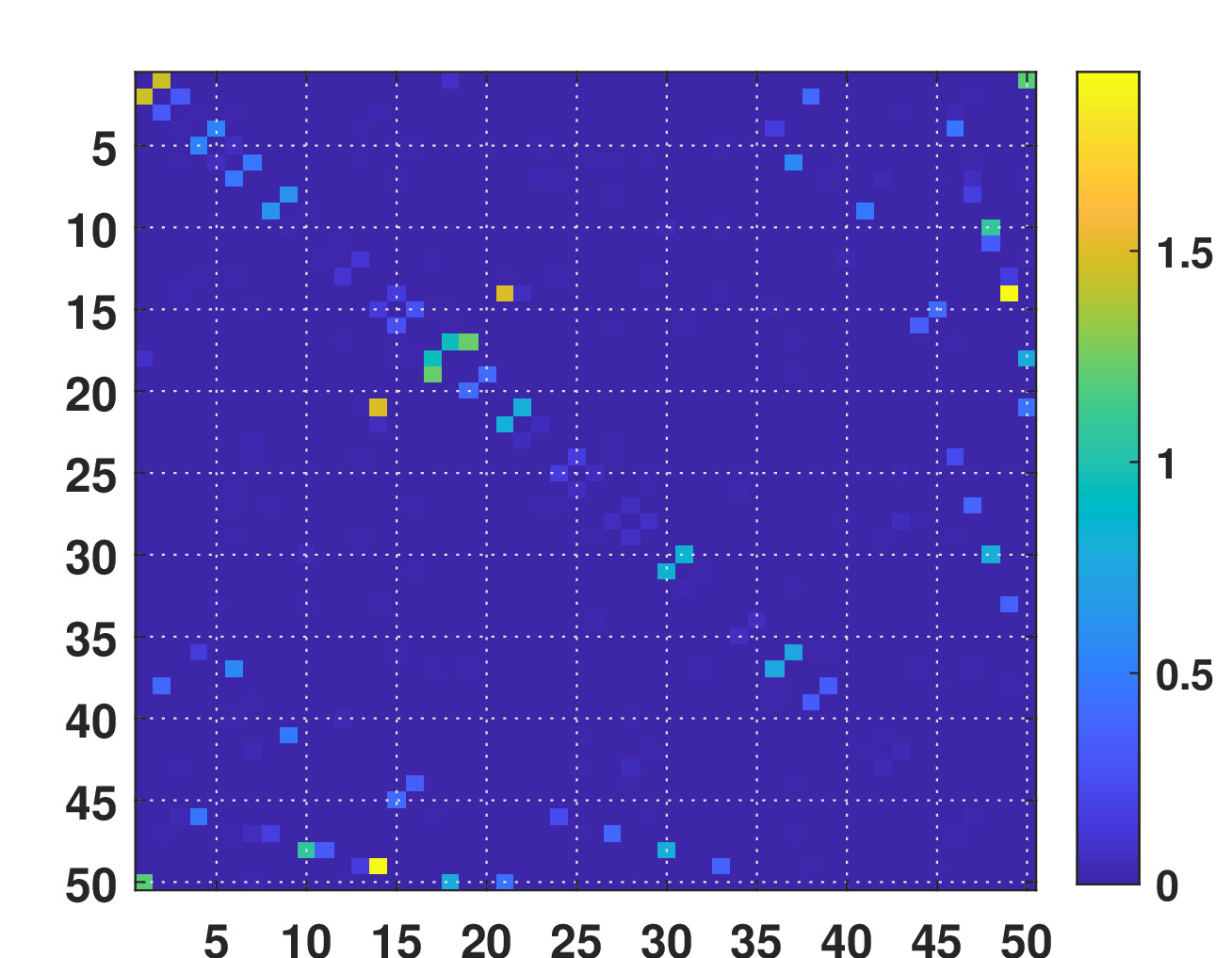} \\
  (a) & (b) \\
	\end{tabular}
	\caption{(a) Generative graph topology $ \mathcal{G_T}(\mathcal{V},  \mathcal{E_T}):$ Under full  observability, green colored edges are included in $ \mathcal{G_T}$. Under partial observability,  red colored nodes are latent, and green colored edges are excluded from $ \mathcal{G_T}.$ (b) Reconstruction under full observability: color map of $H_{\infty}[ \boldsymbol{\hat{\mathsf{W}}}_{ji} ]+H_{\infty}[ \boldsymbol{\hat{\mathsf{W}}}_{ij}]$ for $i,j \in \{1,\cdots,50 \}.$ Some of the spurious edges have very low color intensity. All  edges in $ \mathcal{G_T}$ are identified correctly.  
	\label{fig:Results_topology}}
\end{figure}

\begin{figure}[tb]
	\centering
	\begin{tabular}{cc}
	\hspace{-0.5cm}\includegraphics[width=0.75\columnwidth]{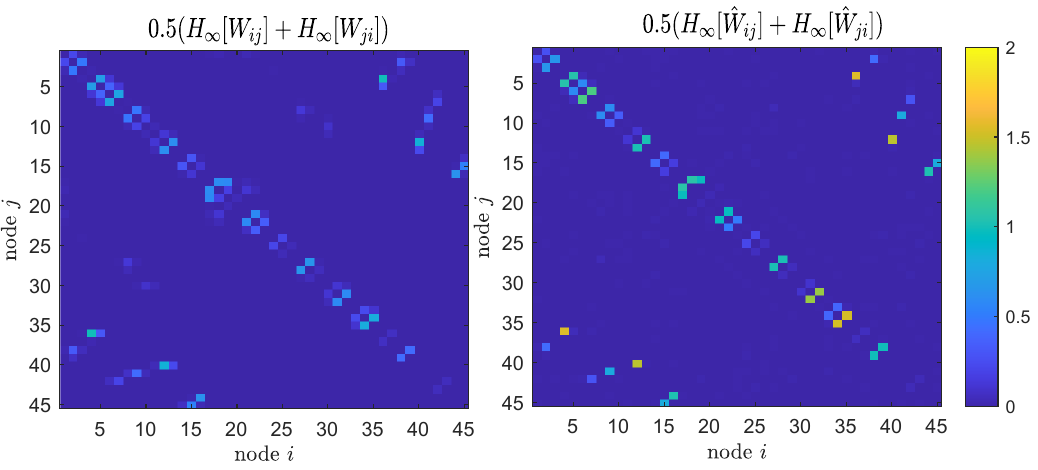} 
	\end{tabular} 
	\caption{Partial observability: Left figure is the magnitude of model-based Wiener filter (exact values based on generative model). Right figure is the magnitude of regularized Wiener filter obtained by solving (\ref{eqn:vec_opt_sparse}) with $\gamma=0.07$ and $L_*=3$. 
	\label{fig:Results_Hnorm_Partial}}
\end{figure}

\begin{figure}[tb]
	\centering
	\begin{tabular}{c}
	\includegraphics[width=0.5\columnwidth]{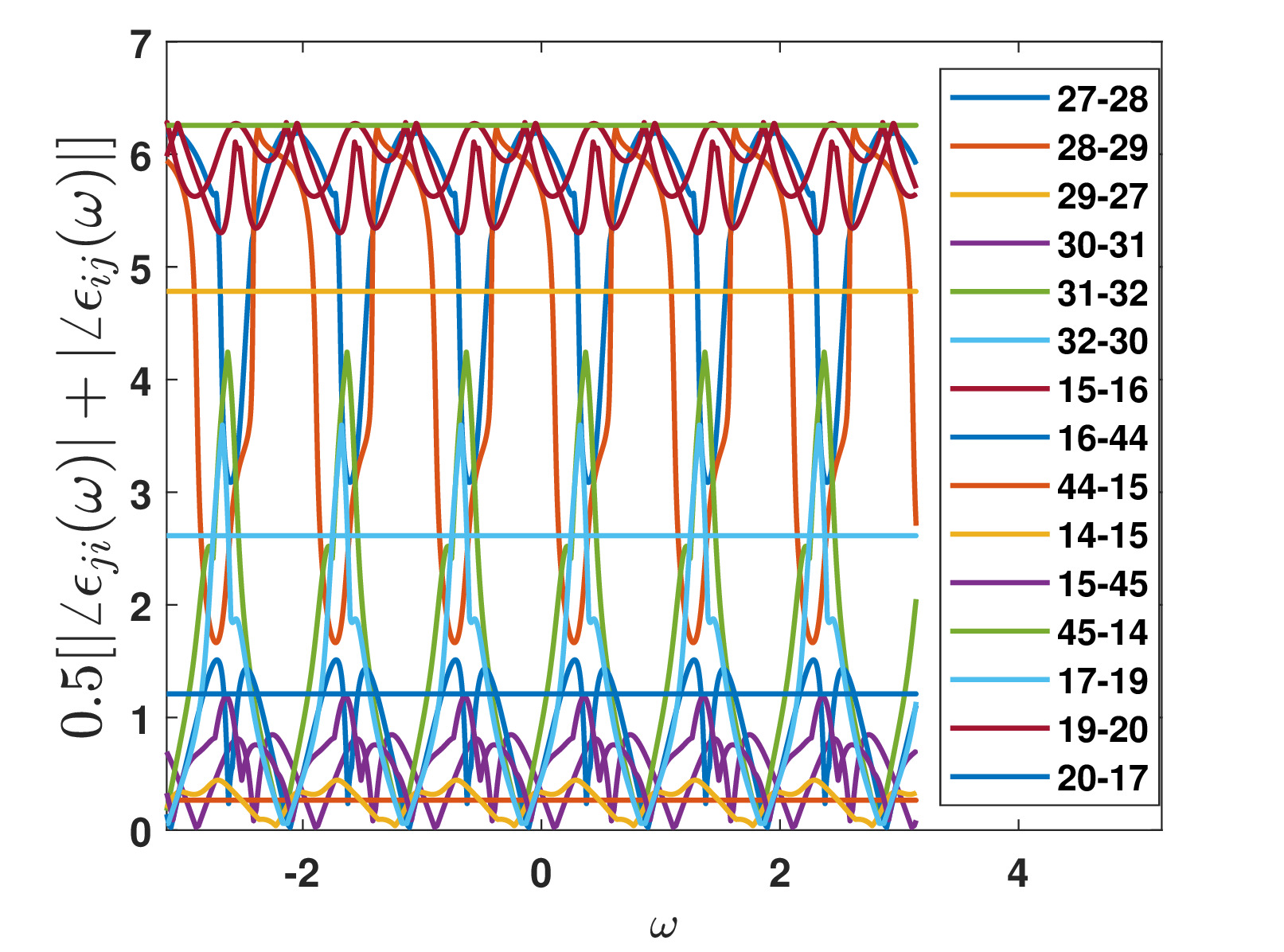}
	\end{tabular}
 \vspace{-0.1in}
	\caption{ Phase plot of eigenvalues (model-based) for various edges $j-i.$ Here, $\epsilon_{ji}$ ( $\epsilon_{ij}$) is the largest component of the vector $\Eps_{ji}$ ($\Eps_{ij}$).
	\label{fig:Results_phase}}
\end{figure}

\begin{figure}[tb]
	\centering
	\includegraphics[width=1\columnwidth]{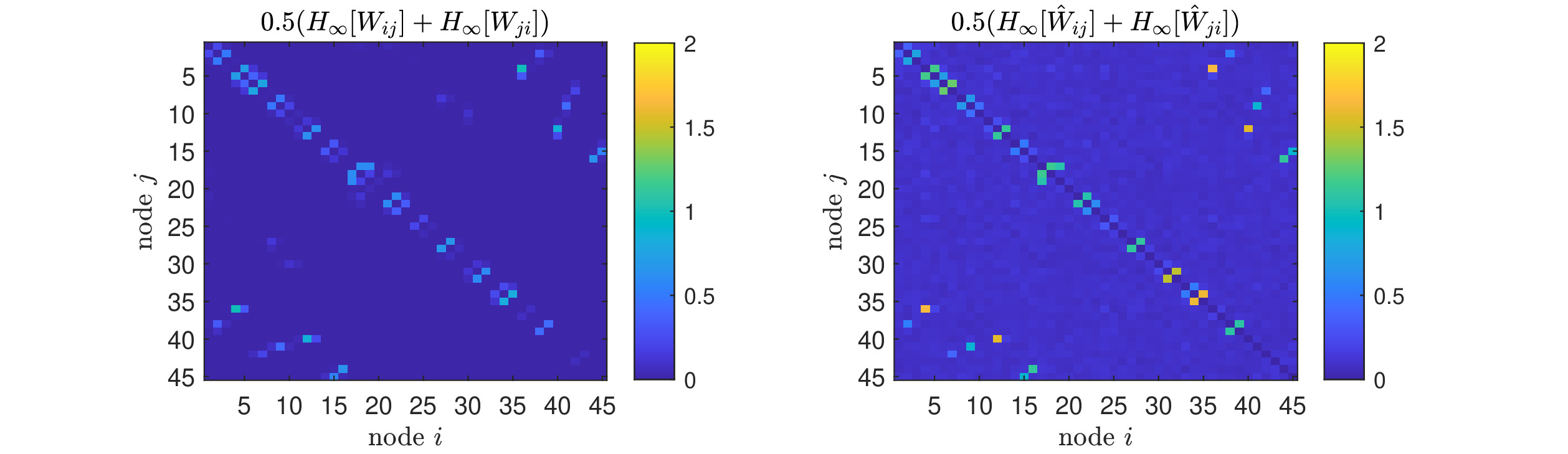} 
	\caption{Reconstruction without regularization: Left figure is the magnitude of model-based Wiener filter (exact values based on generative model). Right figure is the magnitude of Wiener filter estimate without regularizer ($\gamma=0$) and $L_*=3$, computed from 628400 samples.   
	\label{fig:Results_Hnorm_Partial_asymptotic}}
\end{figure}

Please see Fig. \ref{fig:Results_Hnorm_Partial_asymptotic} for consistent topology estimation without regularizer.

\section{Conclusions} \label{sec:six}
For a network of multiple agents interacting according to a linear dynamic model excited by cyclostationary processes, an algorithm is presented for reconstructing the topology using nodal time-series data. 
% The algorithm is based on lifting the time-series to a higher dimension and analyzing the algebraic properties of inverse power spectral density. 
The learning algorithm provably  recovers the true topology in the large sample limit. For bidirected networks with {tree }topologies under partial observability, an algorithm for topology learning is developed. 
% The algorithm is partly based on principle of graphical separation of undirected graphs and phase result. 
The algorithm's performance is demonstrated using test networks, as well as real-time identification of ocean current's  based on sea surface temperature data (COBE-SST2). 
% In future work, we will analyze the sample complexity necessary for consistent estimation as well as understand the effect of measurement noise on the learning processes.

% \section{Epilogue}
% A word or two to conclude, and this even includes some inline 
% maths:  $R(x,t)\sim t^{-\beta}g(x/t^\alpha)\exp(-|x|/t^\alpha)$.

% \begin{ack}                               % Place acknowledgements
% Partially supported by the Roman Senate.  % here.
% \end{ack}

\onecolumn
\appendix

\section{Lemma for \ref{eqn:Z_WSS}:} \label{app:lemma4_proof}
\begin{lem}\label{lem:lemma4}
Consider a collection of cyclostationary processes $\{x_i(k)\}_{i=1}^m$ described by (\ref{eqn:cyclo_LDM}), which are jointly WSCS of period $T$. Consider a lifted process $X_i(k) = [x_i(kT), \cdots, x_i(kT+T-1)]^{\prime}.$ The dynamics of $\{X_i(k)\}_{i=1}^m$ is governed by:
\vspace{-0.5cm}

\small
\begin{align}
X_j(k) &= \sum_{i=1,i \neq j}^{m}(H_{ji}*X_i)(k)+E_j(k).\\ 
\text{The }&z-\text{transform is given by,}\nonumber \\
\mathsf{X}_{j}(z) &= \sum_{i=1,i \neq j}^{m}\mathsf{H}_{ji}(z)\mathsf{X}_i(z) + \mathsf{E}_j(z)\text{; where, }     
\end{align} \normalsize
\small \begin{align} 
E_i(k) &= \begin{bmatrix}e_i(kT)& e_i(kT+1)& \cdots & e_i(kT+T-1) \end{bmatrix}^{\prime},\nonumber\\
X_i(k) &= \begin{bmatrix}x_i(kT) & x_i(kT+1) & \cdots & x_i(kT+T-1) \end{bmatrix}^{\prime},\nonumber\\
\mathsf{H}_{ji}(z) &= D(z^{\frac{1}{T}}) \mathsf{h}_{ji}(z^{\frac{1}{T}}),\ 
D(z) = \begin{bmatrix}z^0&z^1\dots&z^{T-1}\\\vdots&\ddots&\vdots\\z^{-(T-1)}&\dots&z^0\\\end{bmatrix}.\nonumber 
\end{align}
\normalsize
\end{lem}

\begin{pf}

The dynamics of a network of cyclostationary processes are given by,
\small 
\begin{align*}
\mathsf{x}_i(z) &= \sum_{j=1}^{m}\mathsf{h_{ij}}(z)\mathsf{x}_j(z) + \mathsf{e}_i(z),\\
{x}_i(k) &= \sum_{j=1}^{m}\sum_{n= -\infty}^{\infty}h_{ij}(n){x}_j(k-n) + {e}_i(k),\\
{x}_i(kT) &= \sum_{j=1}^{m}\sum_{n= -\infty}^{\infty}h_{ij}(n){x}_j(kT-n) + {e}_i(kT),\\
{x}_i(kT+p) &= \sum_{j=1}^{m}\sum_{n= -\infty}^{\infty}h_{ij}(n){x}_j(kT+p-n) + {e}_i(kT+p),\\
\end{align*} \normalsize

where $k$ is the time index and $p\in \{0,1,\dots,T-1\}$. Substitute $n=aT+b$, where $a\in \mathbb{Z}$ and $b$ takes the values $\{p,p-1,p-2,\dots,p-T+1\}.$
\small \begin{align*}
{x}_i(kT+p) &= \\
\sum_{j=1}^{m}\sum_{n= -\infty}^{\infty}& h_{ij}(n){x}_j(kT+p-n) + {e}_i(kT+p),\\
{x}_i(kT+p) &= \\
\sum_{j=1}^{m}\sum_{a= -\infty}^{\infty}\sum_{b=p-T+1}^{p}& h_{ij}(aT+b){x}_j(kT+p-[aT+b]) + {e}_i(kT+p),\\
{x}_i(kT+p) &=\\ \sum_{j=1}^{m}\sum_{b=p-T+1}^{p}\sum_{a= -\infty}^{\infty}& h_{ij}(aT+b){x}_j(kT+p-[aT+b]) + {e}_i(kT+p),\\
{x}_i(kT+p) &=\\ \sum_{j=1}^{m}\sum_{b=p-T+1}^{p}\sum_{a= -\infty}^{\infty}& h_{ij}(aT+b){x}_j(kT+p-b-aT) + {e}_i(kT+p).\\
\text{Replace }t = p-b\\
{x}_i(kT+p)&=\\
\sum_{j=1}^{m}\sum_{t=0}^{T-1}\sum_{a= -\infty}^{\infty}& h_{ij}(aT+p-t){x}_j([k-a]T+t) + {e}_i(kT+p)\\
\end{align*} \normalsize
Define $H_{ij,pt}[a] =h_{ij}(aT+p-t)$. Here, $p$ is the row index ranging from $0$ to $T-1$ and $t$ is the column index ranging from $0$ to $T-1$ of the filter matrix $H_{ij}$. For the diagonal entries ($p=t$) of $H_{ij}$, $H_{ij,pp}[n] =h_{ij}(nT)$ for $p\in\{0,1,\dots,T-1\}$.
\small \begin{align*}
{x}_i(kT+p) &=\sum_{j=1}^{m}\sum_{t=0}^{T-1}\sum_{a= -\infty}^{\infty} H_{ij,pt}[a]{x}_j([k-a]T+t) + {e}_i(kT+p),\\
{x}_i(kT+p) &=\sum_{j=1}^{m}[ H_{ij,p0}*{x}_j(kT)+\dots+H_{ij,p\ T-1}*{x}_j(kT+T-1)] \\
&\ \ \ +{e}_i(kT+p).
\end{align*} \normalsize

Lift the scalar process $x_i(k)$ to a vector process $X_i(k)$ by varying $p$ from $0$ to $T-1:$ 
\small \begin{align*}
 X_i(k) = \begin{bmatrix} x_i(kT)\\x_i(kT+1)\\\vdots\\x_i(kT+T-1) \end{bmatrix}
\end{align*} \normalsize
Iterate the $p$ from $0$ to $T-1$ to get the following relation:

\small \begin{align*}
 \begin{bmatrix} x_i(kT)\\x_i(kT+1)\\\vdots\\x_i(kT+T-1) \end{bmatrix}=&
 \sum_{j=1}^{m} \begin{bmatrix} H_{ij,00} & H_{ij,01}\ \dots \\ \vdots & \ddots & \\ H_{ij,T-1,0}& & H_{ij,T-1,T-1}
 \end{bmatrix}* \begin{bmatrix} x_j(kT)\\x_j(kT+1)\\\vdots\\x_j(kT+T-1) \end{bmatrix}+\begin{bmatrix} e_i(kT)\\e_i(kT+1)\\\vdots\\e_i(kT+T-1) \end{bmatrix}
\end{align*} 
\begin{align}
 X_i(k) =& \sum_{j=1}^{m}H_{ij}*X_j(k)+E_i(k)\\
 \text{Taking }&\text{the Z transform, we get,} \nonumber \\
 \mathsf{X}_{j}(z) =& \sum_{i=1}^{m}\mathsf{H}_{ji}(z)\mathsf{X}_i(z) + \mathsf{E}_j(z)
\end{align}  \normalsize
where $H_{ij,pt}*x_j(kT+t) = \sum_{a=-\infty}^{\infty} H_{ij,pt}[a]x_j(kT+t-aT)$ for $t\in\{0,1,\dots,T-1\}$ with $H_{ij,pt}[n] =h_{ij}(nT+p-t)=h_{ij}(T[n+\frac{p-t}{T}])$. This implies $\mathsf{H}_{ij}(z) = \mathcal{Z}[H_{ij}(k)] = D(z^{\frac{1}{T}}) \mathsf{h}_{ij}(z^{\frac{1}{T}})$, where $D(z) = \begin{bmatrix}z^0&z^1\dots&z^{T-1}\\\vdots&\ddots&\vdots\\z^{-(T-1)}&\dots&z^0\\\end{bmatrix}$.
\end{pf}

\section{Illustration}\label{app:illustration}

\begin{exmp}
 Consider the following generative model for the time-series, 
\small \begin{align}\label{eqn:lindyn}
 \sum_{n=0}^{l}a_{n,i}\frac{d^nx_i}{dt^n}=\sum_{j=1,j\neq i}^{m}b_{ij}x_j(t) + p_i(t), 
\end{align} \normalsize
where, $a_{n,i}, \ b_{ij}$ are complex-valued system constants, $p_i(t)$ is an exogenous input, $x_i(t)$ is the state of $i^{th}$ agent. Representing (\ref{eqn:lindyn}) in the form of (\ref{eqn:cyclo_LDM}) using bilinear transform (Tustin's method \cite{oppenheim1999discrete}), we get 
$\mathsf{h}_{ij}(z) = \frac{b_{ij}}{\mathsf{S_i}(z)}$, $\mathsf{e}_i(z):=$ $\frac{\mathsf{p}_i(z)}{\mathsf{S_i(z)}},$ 
where $\Delta t$ is the sampling period and $\mathsf{S_i}(z):=\sum_{n=1}^{l}a_{n,i}(\frac{2(1-z^{-1})}{\Delta t(1+z^{-1})})^{n}.$ {{Assumption \ref{ass:transferFunction_property} holds for this example with $a_{jki} = \frac{b_{kj}}{b_{ki}}$ for all $i,j,k$ such that $b_{kj},b_{ki} \neq 0$ holds in (\ref{eqn:lindyn}).}} 
\end{exmp}

\section{Pathological cases:} \label{app:thm3.3_proof}

\begin{thm}\label{thm:3.3}
Consider a well-posed and topologically detectable networked LDS described by (\ref{eqn:cyclo_LDM}), with its equivalent representation $(\mathbb{H}(z),{E})$ as in (\ref{vector_LDM_gen}), with associated graph $\mathcal{G},$ topology $\mathcal{G_T}$ and satisfying Assumptions \ref{ass:transferFunction_property} and \ref{ass:atmost1}. Suppose $i$ and $j$ are such that $i \in \mathcal{N}_{\mathcal{G_T}}(j)$, and denote the eigenvalues of ${\boldsymbol{B}_{j}^{\prime}}{{\Phi}_{X}^{-1}}{\boldsymbol{B}_{i}}$ as $\Eps_{ji}\in \C^{T\times 1}$. 
The set of system parameters $\{a_{jki}, \forall k \in \mathcal{V}\setminus i,j\},$ such that $\phase{\Eps_{ji}(l)}$ is a constant for $l\in \{1,2,\cdots,T \}$ and all $\omega \in [0, 2 \pi),$ has measure zero.
\end{thm}

\begin{pf} Using (\ref{PhiInv_expression}) and proof of Theorem \ref{thm:3.2}, we have
\small \begin{align*}
&{\boldsymbol{B}_{j}^{\prime}}{{\Phi}_{X}^{-1}(\omega)}{\boldsymbol{B}_{i}}\\
&= -b_{ji}\mathsf{S^{*}_j}(\frac{\omega}{T}){{\Phi}_{P_j}^{-1}(\omega)}D(\frac{\omega}{T})-b_{ij}\mathsf{S_i}(\frac{\omega}{T})D^{*}(\frac{\omega}{T}){{\Phi}_{P_i}^{-1}(\omega)}\\
&\ \ \ \ \ \ \ \ \ \ \ \ \ \ \ +{(b_{kj}b_{ki})}[D^*(\frac{\omega}{T}){\Phi}_{E_k}^{-1}(\omega)D(\frac{\omega}{T})].\\
\end{align*} \normalsize

Suppose $A,B$ and $C$ correspond to first, second and third term in the above equation respectively. Given that the phase response of the eigenvalues of ${\boldsymbol{B}_{j}^{\prime}}{{\Phi}_{X}^{-1}(\omega)}{\boldsymbol{B}_{i}}$ is a constant and is denoted by $c_{ji}.$ The eigenvalues of $C$ has a constant phase of $\phase{b_{kj}b_{ki}}.$ Thus, the phase of the eigenvalues of $A+B$ is another constant $\alpha = c_{ji}- \phase{b_{ki}b_{kj}}.$ Thus, $e^{-\iota \alpha }[A+B]$ has to be Hermitian and has real eigenvalues. It follows from the property of Hermitian that,
\small
\begin{align*}
&e^{-\iota \alpha}[-b_{ji}\mathsf{S^{*}_j}(\frac{\omega}{T}){{\Phi}_{P_j}^{-1}(\omega)}D(\frac{\omega}{T})-b_{ij}\mathsf{S_i}(\frac{\omega}{T})D^{*}(\frac{\omega}{T}){{\Phi}_{P_i}^{-1}(\omega)}]\\
&=e^{\iota \alpha}[-b^{*}_{ji}\mathsf{S_j}(\frac{\omega}{T})D^{*}(\frac{\omega}{T}){{\Phi}_{P_j}^{-1}(\omega)}-b^{*}_{ij}\mathsf{S^{*}_i}(\frac{\omega}{T}){{\Phi}_{P_i}^{-1}(\omega)}D(\frac{\omega}{T})].\\
\end{align*}
\normalsize

The set of system parameters which satisfies the above equality constraint for all $\omega \in [0,2\pi)$ has a zero measure. 
\end{pf}

\section{Algorithm $1$} \label{app:alg1}

\small
\begin{algorithm}
\caption{Learning Algorithm for reconstructing the topology of LDG with cyclostationary inputs}
\textbf{Input:} Nodal time-series $x_i(k)$ for each node $i \in \{1, 2,... m\}$ which is WSCS. Thresholds $\rho,\tau$. Frequency points $\Omega$. \\
\textbf{Output:} Reconstructed Topology ($\mathcal{V},\bar{\mathcal {E_T}}$) \\
\begin{algorithmic}[1]
\State Compute the periodogram of $x_i(k)$ to determine the period $T_i$. Determine $T= LCM\{T_1,\cdots, T_m \}.$ Lift each $x_i(k)$ to vector process $X_i(k)$ of size $T.$ Define $X(k) = [X_1(k),\cdots X_m(k)]^{\prime}$
\State Edge set $\bar{\mathcal{E}}_M \gets \{\}$ 
\ForAll{$l,p \in \{1,2,...,m\}, l\neq p$}
\If{$H_{\infty}[ \boldsymbol{\mathsf{W}}_{lp} ] \neq 0$}
\State $\bar{\mathcal{E}}_M \gets \bar{\mathcal{E}}_M \cup \{(l,p)\}$
\EndIf
\EndFor\label{step1_b}
\State Edge set $\bar{\mathcal{E}}_T \gets \bar{\mathcal{E}}_M$ 
\ForAll{$(p,l) \in \bar{\mathcal{E}}_M$}
\State Compute $\{{\Eps_{pl}(t)\}}^{\prime}_{t=1}=eig\{\boldsymbol{\mathsf{W}}_{pl}(\omega)\}$
\If{$\phase{\Eps_{pl}(t)}$ is constant $\forall \omega \in [0,2\pi), \forall t$}
\State $\bar{\mathcal{E}}_T \gets \bar{\mathcal{E}}_T - \{(p,l)\}$
\EndIf
\EndFor 
\end{algorithmic}
\end{algorithm}
\normalsize

\newpage
\section{Algorithm $2$} \label{app:alg2}

\small
\begin{algorithm}
\small{
\caption{Learning $\mathcal{G}_{T_o}(\mathcal{V}_o,\mathcal{E}_{T_o})$}
\textbf{Input:} $\mathcal{G}_c=(\mathcal{V}_o,\mathcal{E}_c)$ generated from structure of $\Phi_{oo}^{-1}$\\
\textbf{Output:} $\overline{\mathcal{T}}=(\mathcal{V}_o, \mathcal{E}_{\overline{\mathcal{T}}}),$ $V_l$ and $V_{nl}.$
\begin{algorithmic}[1]
\State Edge set ${\mathcal{E}}_{\overline{\mathcal{T}}} \gets \{\}$
\State $V_{nl} \gets \{ \}$
\ForAll{edge $a - b$ in $\mathcal{E}_c$}
\If{$Z:=\{a,b\}$ there exist $I \neq \{\phi\}$ and $J \neq \{\phi\}$ such that $sep(I,Z,J)$ holds in $\mathcal{G}_{c}$}
\State $V_{nl} \gets V_{nl} \cup \{a,b\}$, {$\mathcal{E}_{\overline{\mathcal{T}}} \gets \mathcal{E}_{\overline{\mathcal{T}}} \cup \{(a,b)\}$
\EndIf
\EndFor
\State $V_l \gets \mathcal{V}_o - V_{nl}$
%\ForAll{$a \in V_{l}, b \in V_{nl}$ with $(a,b) \in \mathcal{E}_{\mathcal{G}_c}$}
%\State $\mathcal{E}_{\overline{\mathcal{T}}} \gets \mathcal{E}_{\overline{\mathcal{T}}} \cup \{(a,b)\}$
%\EndFor
\ForAll{$a \in V_{l}, b \in V_{nl}$ with $(a,b) \in \mathcal{E}_{\mathcal{G}_c}$}
\State Compute $\Eps_{ab}= eig[\Phi^{-1}_{oo}(a,b)]$
\If{$\phase{\Eps_{ab}[t]} $ is not a constant $\forall \omega \in [0,2\pi)$, for $t \in \{1,\cdots,T\}$}
\State $\mathcal{E}_{\overline{\mathcal{T}}} \gets \mathcal{E}_{\overline{\mathcal{T}}} \cup \{(a,b)\}$
\EndIf
\EndFor\label{step1_c}}
\end{algorithmic}}
\end{algorithm}
\normalsize

\section{Algorithm $3$} \label{app:alg3}
\small
\begin{algorithm}
\small{
\caption{Reconstructing $\mathcal{V}_h$ and $\mathcal{E}_{T_h}$}
\textbf{Input:} $\overline{\mathcal{T}}=(\mathcal{V}_o,\mathcal{E}_{\overline{\mathcal{T}}})=\cup_{j=1}^{h}\overline{\mathcal{T}}_j$\\
\textbf{Output:} $\tilde{\mathcal{T}} = (V_{\tilde{\mathcal{T}}},\mathcal{E}_{\tilde{\mathcal{T}}})$.
\begin{algorithmic}[1]
\State Node set $V_{\tilde{\mathcal{T}}} \gets \mathcal{V}_o$, edge set $\mathcal{E}_{\tilde{\mathcal{T}}}\gets \mathcal{E}_{\overline{\mathcal{T}}}$
\State $h\gets$ Number of disjoint subgraphs in $\mathcal{G}_{T_o}$
\ForAll{$j \in \{1,2,...,h\}$}
\ForAll{$i \in \{j+1,...,h\}$}
\If{there exist a pair of nodes $a,b$ such that $a \in \overline{\mathcal{T}}_j$ and $b \in \overline{\mathcal{T}}_i$ such that all their neighbors in $\overline{\mathcal{T}}$ are connected in $\mathcal{G}_c$}
\State $V_{\tilde{\mathcal{T}}} \gets V_{\tilde{\mathcal{T}}} \cup l_j$
\State $\mathcal{E}_{\tilde{\mathcal{T}}} \gets \mathcal{E}_{\tilde{\mathcal{T}}} \cup \{(a,l_j),(l_j,b)\}$
\EndIf
\EndFor\label{step1_b}
\EndFor
\State Merge hidden nodes that are neighbors of the same observed node.
\end{algorithmic}}
\end{algorithm}
\normalsize

\end{document}